# A Nonlinear Heat Equation Arising from Automated-Vehicle Traffic Flow Models

**Dionysios Theodosis**[*], **Iasson Karafyllis**[**], **George Titakis**[*], **Ioannis Papamichail**[*] **and Markos Papageorgiou**[*,***]

[*] Dynamic Systems and Simulation Laboratory,
Technical University of Crete, Chania, 73100, Greece,

emails: dtheodosis@dssl.tuc.gr, gtitakis@dssl.tuc.gr,
ipapa@dssl.tuc.gr, markos@dssl.tuc.gr

[**] Dept. of Mathematics, National Technical University of Athens,
Zografou Campus, 15780, Athens, Greece,
emails: iasonkar@central.ntua.gr , iasonkaraf@gmail.com

[***] Faculty of Maritime and Transportation,
Ningbo University, Ningbo, China.

## Abstract

In this paper, a new nonlinear heat equation is studied that arises as a model of the collective behavior of automated vehicles. The properties of the solutions of this equation are studied by introducing the appropriate notion of a weak solution that requires certain entropy-like conditions. To obtain an approximation of the solution of the nonlinear heat equation, a new conservative first-order finite difference scheme is proposed that respects the corresponding entropy conditions, and certain links between the weak solution and the numerical scheme are provided. Finally, a traffic simulation scenario and a comparison with the Lighthill-Witham-Richards (LWR) model are provided, illustrating the benefits of the use of automated vehicles.

## 1. Introduction

The heat equation is a fundamental differential equation and among the most studied topics, see for instance [11], [20], [21], [25], [26], [15], [27], [30], [33] and references therein. It describes the evolution of the density of some quantity, such as temperature or chemical concentration, and finds immediate application in a variety of topics such as financial mathematics, flow and pressure diffusion in porous media, image analysis and many others.

In this paper, we consider a novel nonlinear heat equation, whose diffusion/viscosity coefficient depends on both the density and the spatial derivative of the density (see also [11], [20], [21], [15]), a feature that is rarely studied in the literature where typically the diffusion coefficient exclusively depends on the density, see for instance [11], [25], [26], [27], [29], [30], [33] and references therein. The nonlinear heat equation under consideration is inspired by a second order macroscopic traffic flow model that was recently derived in [16] by the design of cruise controllers for autonomous vehicles and has several fluid-like characteristics. Macroscopic traffic flow models have been widely used as they are able to capture the collective behavior of vehicles in a traffic stream and can express certain relationships between traffic flow, density and mean speed. Macroscopic traffic flow models can be in general distinguished between first-order models, governed by the continuity equation, and second-order models, which include an additional differential equation for the speed (see for instance [2], [12], [34]).

In the present work, it is first shown that the solutions of the corresponding second-order macroscopic model can be approximated by the solutions of this particular nonlinear heat equation. Furthermore, inspired by the mechanical energy of the original second-order system, certain functionals are defined, providing entropy-like conditions that characterize physically meaningful solutions, see [9], [31]. In particular, we define the concept of a weak solution for the problem at hand that is more demanding than the typical definition of weak solution (see [24], [31]), as it requires that certain entropy conditions hold and, in addition, the mass is conserved. Finally, it is shown that, under certain regularity assumptions, this type of weak solutions is equivalent to classical solutions, under the assumption of certain entropy-like conditions (see Proposition 1).

The study of the nonlinear heat equation is inevitably related to numerical techniques in order to obtain an approximation of its solution. More specifically, several finite-difference methods have been proposed to study a variety of nonlinear heat equations, see for instance [3], [4], [6], [7], [8], [10], [13], [14], and references therein. In this paper, we propose a new first-order explicit finite-difference scheme that, in addition to being conservative (Proposition 2), it also respects the corresponding entropy conditions (Proposition 3 and Proposition 5). Moreover, a condition on the step size is provided that is a nonlinear version of the Courant-Friedrichs-Levy condition (Proposition 4). Finally, a link between the numerical solution produced by the proposed numerical scheme and the corresponding weak solution is provided (Proposition 6). To demonstrate the properties of the solutions of the nonlinear heat equation, we present two examples: (i) an academic example, where the efficiency of the proposed numerical scheme is illustrated; and (ii) a realistic traffic simulation scenario showing that the mean flow produced by automated vehicles is much higher than that of the LWR (see [22], [28]).

Therefore, the contribution of the paper is summarized as follows:
- The study of a novel nonlinear heat equation that arises in traffic flow theory for automated vehicles;
- The construction of a first-order numerical scheme that can be applied to nonlinear heat equations with an explicit condition for numerical stability;
- The study of the connection of the weak solution for the nonlinear heat equation with the solution produced by the proposed numerical scheme.

The structure of the paper is as follows. Section 2 is devoted to the presentation of a second-order traffic flow model and the approximation of its solution by the solution of a particular nonlinear heat equation. Section 3 presents the proposed numerical scheme and its properties. Section 4 presents the two aforementioned examples. Finally, some concluding remarks are given in Section 5. All proofs are provided in the Appendix.

**Notation.** Throughout this paper, we adopt the following notation.

* By $|x|$ we denote both the Euclidean norm of a vector $x \in \Re^n$ and the absolute value of a scalar $x \in \Re$.

* $\Re_+ := [0, +\infty)$ denotes the set of non-negative real numbers.

* By $x'$ we denote the transpose of a vector $x \in \Re^n$. By $|x|_\infty = \max\{|x_i|, i = 1, ..., n\}$ we denote the infinity norm of a vector $x = (x_1, x_2, ..., x_n)' \in \Re^n$.

* Let $A \subseteq \Re^n$ be an open set. By $C^0(A; \Omega)$, we denote the class of continuous functions on $A \subseteq \Re^n$, which take values in $\Omega \subseteq \Re^m$. By $C^k(A; \Omega)$, where $k \geq 1$ is an integer, we denote the class of functions on $A \subseteq \Re^n$ with continuous derivatives of order $k$, which take values in $\Omega \subseteq \Re^m$. When $\Omega = \Re$ the we write $C^0(A)$ or $C^k(A)$.

* Let $I \subseteq \Re$ be a given interval. For $p \in [1, \infty)$, $L^p(I)$ denotes the set of equivalence classes of Lebesgue measurable functions $f : I \to \Re$ with $\|f\|_p := \left( \int_I |f(x)|^p \, dx \right)^{1/p} < +\infty$. $L^\infty(I)$ denotes the set of equivalence classes of measurable functions $f : I \to \Re$ for which $\|f\|_\infty = ess \sup_{x \in I} (|f(x)|) < +\infty$.

* Let $u : \Re_+ \times \Re \to \Re$, $(t, x) \to u(t, x)$ be any function differentiable with respect to its arguments. We use the notation $u_t(t, x) = \frac{\partial u}{\partial t}(t, x)$ and $u_x(t, x) = \frac{\partial u}{\partial x}(t, x)$ for the partial derivatives of $u$ with respect to $t$ and $x$, respectively. We use the notation $u[t]$ to denote the profile at certain $t \geq 0$, $(u[t])[x] := u(t, x)$, for all $x \in \Re$.

* For a set $S \subseteq \Re^n$, $\bar{S}$ denotes the closure of $S$.

## 2. A Nonlinear Heat Equation

*2.1 Motivation and Derivation of a Nonlinear Heat Equation*

In the recent paper [16], the study of microscopic vehicle movement control laws (cruise controllers) for autonomous vehicles on lane-free roads brought forth the following macroscopic model that holds for $\tau > 0$, $\xi \in \Re$:

$$\tilde{\rho}_\tau + (\tilde{\rho}\tilde{v})_\xi = 0 \tag{2.1}$$

$$\tilde{\rho}\tilde{q}(\tilde{v})\tilde{v}_\tau + \tilde{\rho}\tilde{q}(\tilde{v})\tilde{v}\tilde{v}_\xi + \tilde{P}'(\tilde{\rho})\tilde{\rho}_\xi = (\tilde{\mu}(\tilde{\rho})\tilde{g}'(\tilde{v})\tilde{v}_\xi)_\xi - \tilde{\rho}\tilde{f}(\tilde{v} - v^*) \tag{2.2}$$

with constraints $\tilde{\rho}(\tau, \xi) \in [0, \rho_{max})$, $\tilde{v}(\tau, \xi) \in (0, v_{max})$ for $\tau > 0$, $\xi \in \Re$, where $\rho_{max} > 0$, $v_{max} > 0$, and $v^* \in (0, v_{max})$ are constants. The states $\tilde{\rho}(\tau, \xi)$, $\tilde{v}(\tau, \xi)$ are the traffic density and mean speed, respectively, at time $\tau > 0$ and position $\xi \in \Re$ on a highway, while the constants $\rho_{max} > 0$, $v_{max} > 0$, and $v^* \in (0, v_{max})$ are the maximum density, the maximum speed, and the speed set-point, respectively. Moreover, $\tilde{f} : \Re \to \Re$ is a $C^1$ function with $\tilde{f}(0) = 0$ and $\xi \tilde{f}(\xi) > 0$ for all $\xi \neq 0$, $\tilde{g} \in C^1(\Re)$ is an increasing function with $\tilde{g}'(\tilde{v}) > 0$ for all $\tilde{v} \in \Re$ and $\tilde{q}(\tilde{v})$ is defined by

$$\tilde{q}(\tilde{v}) = v_{max}^2 \frac{v_{max}\tilde{v} - 2\tilde{v}v^* + v^* v_{max}}{2(v_{max} - \tilde{v})^2 \tilde{v}^2} \tag{2.3}$$

Model (2.1), (2.2) is a fluid-like model where the term $\tilde{P}'(\tilde{\rho})\tilde{\rho}_\xi$ is a pressure term and expresses the tendency to accelerate or to decelerate based on the (local) density; while the term $(\tilde{\mu}(\tilde{\rho})\tilde{g}'(\tilde{v})\tilde{v}_\xi)_\xi$, is a viscosity term, by analogy with the theory of fluids, with $\tilde{\mu}(\tilde{\rho})$ playing the role of dynamic viscosity, see [17], [23]. The functions $\tilde{\mu} : [0, \rho_{max}) \to \Re_+$, $\tilde{P} : [0, \rho_{max}) \to \Re_+$ are $C^1([0, \rho_{max}))$ and satisfy the following properties

$$\lim_{\tilde{\rho} \to \rho_{max}^-} \tilde{P}(\tilde{\rho}) = +\infty \tag{2.4}$$

$\tilde{\mu}(\tilde{\rho}) = 0$, $\tilde{P}(\tilde{\rho}) = 0$ for all $\tilde{\rho} \in [0, \bar{\rho}]$, and $\tilde{\mu}(\tilde{\rho}) > 0$, $\tilde{P}(\tilde{\rho}) > 0$ for all $\tilde{\rho} \in (\bar{\rho}, \rho_{max})$ (2.5)

where the constant $\bar{\rho} \in (0, \rho_{max})$ is referred to as the interaction density due to the fact that for all $\rho \in [0, \bar{\rho}]$ there is no interaction among the vehicles (see (2.5)). The term $f(\tilde{v} - v^*)$ is a *relaxation term* that describes the tendency of vehicles to adjust their speed to the given speed set-point $v^* \in (0, v_{max})$. The functions $\tilde{\mu}(\tilde{\rho})$, $\tilde{g}(\tilde{\rho})$, $f(\tilde{v} - v^*)$, $\tilde{P}(\tilde{\rho})$ are determined by the properties of the cruise controller of the vehicles (see [16]).

**Remarks: (i)** Model (2.1), (2.2) does not include non-local terms and has certain characteristics from the kinematic theory of fluids. Traffic flow is isotropic, as in fluid flow, since the vehicles are autonomous and do not react to downstream vehicles only (as in conventional traffic), but also to upstream vehicles.
**(ii)** There are infinite equilibrium points, namely the points where $\tilde{v}(\xi) \equiv v^*$ and $\tilde{\rho}(\xi) \leq \bar{\rho}$ for all $\xi \in \Re$.
**(iii)** The selection of $\tilde{\mu}(\tilde{\rho})$, $\tilde{g}(\tilde{\rho})$ has several implications for the characteristics of the traffic flow. The dynamic viscosity $\tilde{\mu}(\tilde{\rho})$ makes the "traffic fluid" act as a Newtonian fluid. However, in contrast to actual fluids, the dynamic viscosity also satisfies (2.5). For isentropic (or barotropic) flow of gases, the dynamic viscosity and the pressure are always increasing functions of the fluid density (see the discussion in [17]). Thus, a traffic fluid can have very different physical properties from those of real compressible fluids.

Model (2.1), (2.2) is to be studied under the following conditions for all $\tau \geq 0$:

$$\lim_{\xi \to \pm\infty} \tilde{\rho}(\tau,\xi) = 0 \tag{2.6}$$

$$\lim_{\xi \to \pm\infty} \tilde{v}(\tau,\xi) = v^* \tag{2.7}$$

$$\int_{-\infty}^{+\infty} \tilde{\rho}(\tau,\xi)d\xi < +\infty . \tag{2.8}$$

Condition (2.6) expresses the fact that "far downstream and far upstream" the highway is "empty", i.e., there are no vehicles. Condition (2.8) expresses the fact that the total mass of the vehicles in the highway is finite. Let $r \in \Re$ be a given constant length. Using the variable transformation $\xi = rx + v^*\tau$, $x \in \Re$ and the dimensionless quantities $\tau = \frac{r}{v^*}t$, $b = \frac{v_{max} - v^*}{v^*}$, $R = \frac{\rho_{max}}{\rho} > 1$, $w = \frac{\tilde{v} - v^*}{v^*}$, $\rho = \frac{\tilde{\rho}}{\rho}$, we obtain the following dimensionless model for $t > 0$, $x \in \Re$

$$\rho_t + (\rho w)_x = 0 \tag{2.9}$$

$$\rho q(w) w_t + \rho q(w) w w_x + P'(\rho)\rho_x = (\mu(\rho) g'(w) w_x)_x - \rho f(w) \tag{2.10}$$

With constraints $\rho(t,x) \in [0,R)$, $w(t,x) \in (-1,b)$ for $t > 0$, $x \in \Re$, where $R > 1$, $b > 0$ are constants, $g \in C^1(\Re)$ is an increasing function with $g'(w) > 0$ for all $w \in (-1,b)$, $f: \Re \to \Re$ with $f(0) = 0$ is a $C^1$ function with $wf(w) > 0$ for all $w \neq 0$, and

$$q(w) = (1+b)^2 \frac{2b + (b-1)w}{2(b-w)^2 (1+w)^2}, \text{ for all } w \in (-1,b). \tag{2.11}$$

Moreover, the functions $\mu:[0,R) \to \Re_+$, $P:[0,R) \to \Re_+$ are $C^1([0,R))$ and satisfy the following properties:

$$\lim_{\rho \to R^-} P(\rho) = +\infty \tag{2.12}$$

$\mu(\rho) = 0$, $P(\rho) = 0$ for all $\rho \in [0,1]$ and $\mu(\rho) > 0, P(\rho) > 0$ for all $\rho \in (1,R)$. (2.13)

Model (2.9), (2.10) is studied under the following conditions for all $t \geq 0$:

$$\lim_{x \to \pm\infty} \rho(t,x) = 0 \tag{2.14}$$

$$\lim_{x \to \pm\infty} w(t,x) = 0 \tag{2.15}$$

$$\int_{-\infty}^{+\infty} \rho(t,x)dx < +\infty \tag{2.16}$$

If the parameters of the cruise controller that is applied to the vehicles are selected in such a way that

$$g(w) = w, \text{ for all } w \in (-1,b) \tag{2.17}$$

$$P'(\rho) = k\rho^{-1}\mu(\rho) \tag{2.18}$$

$$\beta(w) = \int_0^w q(s)ds = \frac{b+1}{2}\left[\frac{w(b+1)}{(w+1)(b-w)} + \ln\left(\frac{b(w+1)}{b-w}\right)\right], \text{ for all } w \in (-1,b) \tag{2.19}$$

$$f(w) = k\beta(w) \tag{2.20}$$

where $k > 0$ is a constant and the function $\mu:[0,R) \to \Re_+$ satisfies

$$\lim_{\rho \to R^-} \mu(\rho) = +\infty \tag{2.21}$$

then using definitions (2.17), (2.18), (2.19), and (2.20), the resulting PDE model (2.9), (2.10), becomes

$$\rho_t + (\rho w)_x = 0 \tag{2.22}$$

$$\rho\beta'(w)w_t + \rho\beta'(w)ww_x + k\rho^{-1}\mu(\rho)\rho_x = (\mu(\rho)w_x)_x - k\rho\beta(w) \tag{2.23}$$

with $t > 0$, $x \in \Re$.

**Remark:** Definitions (2.17), (2.20), (2.18), (2.21), in addition to simplifying the original model (2.9), (2.10), have important implications. Definition (2.17) is met in real fluid flows in non-porous media, while condition (2.21) for the dynamic viscosity makes the fluid behave like a solid when density tends to the maximum density $R$.

Assume now that $\rho(t,x) > 0$ for $t > 0$, $x \in \Re$ and that (2.14), (2.15), and (2.16) hold. Define

$$\varphi(t,x) = \rho(t,x)\beta(w(t,x)) + \rho^{-1}(t,x)\mu(\rho(t,x))\rho_x(t,x), \text{ for all } t > 0 \text{ and } x \in \Re \tag{2.24}$$

where $\beta(w)$ is given in (2.19). Equations (2.22), (2.23) and definition (2.24) imply the following equation for $t > 0$, $x \in \Re$:

$$\varphi_t + (w\varphi)_x = -k\varphi. \tag{2.25}$$

Next, define

$$g(t,x) = \frac{\varphi(t,x)}{\rho(t,x)}\exp(kt), \text{ for all } t > 0 \text{ and } x \in \Re. \tag{2.26}$$

Thus, we get from (2.22), (2.23), (2.25), and (2.26) that

$$g_t + w g_x = 0. \tag{2.27}$$

It follows from (2.27) that $\|g[t]\|_\infty \leq \|g[0]\|_\infty$ for all $t \geq 0$, which, combined with (2.26), implies that the following estimate holds:

$$\left\| \frac{\varphi[t]}{\rho[t]} \right\|_\infty \leq \exp(-kt) \left\| \frac{\varphi[0]}{\rho[0]} \right\|_\infty \quad \text{for all } t \geq 0. \tag{2.28}$$

Inequality (2.28) and definition (2.24) imply that for all $x \in \Re$:

$$\lim_{t \to +\infty} \left( \beta(w(t,x)) + \rho^{-2}(t,x) \mu(\rho(t,x)) \rho_x(t,x) \right) = 0. \tag{2.29}$$

Since $\beta((-1,b)) = \Re$, we can define the inverse function $h : \Re \to (-1, b)$ with $h(0) = 0$ and $h'(s) > 0$ for all $s \in \Re$. Thus, we get from (2.29) for all $x \in \Re$:

$$\lim_{t \to +\infty} \left( w(t,x) - h\left( -\rho^{-2}(t,x) \mu(\rho(t,x)) \rho_x(t,x) \right) \right) = 0. \tag{2.30}$$

Inequality (2.28) and equation (2.30) imply that the manifold $w(t,x) = h\left( -\rho^{-2}(t,x) \mu(\rho(t,x)) \rho_x(t,x) \right)$ is exponentially attracting.

Consequently, the solutions of (2.9), (2.10) are approximated by the solutions of the following PDE that holds for $t > 0$, $x \in \Re$

$$\rho_t + \left( \rho h(-\kappa(\rho) \rho_x) \right)_x = 0$$

where

$$\kappa(\rho) := \rho^{-2} \mu(\rho) \quad \text{for } \rho \in (0, R). \tag{2.31}$$

Notice that due to definition (2.31), (2.13) and (2.21), it follows that $\kappa(\rho) = 0$ for all $\rho \in (0, 1]$ and $\kappa(\rho) > 0$ for all $\rho \in (1, R)$.

## 2.2 The Model and its Properties

Motivated by the analysis of the previous section, we study hereafter the following equation and its properties

$$\rho_t + \left( \rho h(-\kappa(\rho) \rho_x) \right)_x = 0, \quad t > 0, \ x \in \Re \tag{2.32}$$

where $\kappa \in C^1([0,R]; \Re_+)$ with $\kappa(\rho) = 0$ for all $\rho \in [0,1]$ and $\kappa(\rho) > 0$ for all $\rho \in (1,R)$, and $h \in C^1(\Re;(-1,b))$ with $h' \in L^\infty(\Re)$ and satisfies $h(0) = 0$ and $h'(s) > 0$ for all $s \in \Re$. In addition, the PDE (2.32) is studied under the condition:

$$\lim_{x \to \pm\infty} (\rho(t,x)) = 0, \text{ for all } t > 0 \tag{2.33}$$

with state constraint $\rho(t,x) \in [0,R)$ for all $t > 0$ and $x \in \Re$, and $\int_{-\infty}^{+\infty} \rho(t,x) dx < +\infty$ for all $t \geq 0$.

We notice that :

- There are infinite equilibrium points, namely the points where $\rho(x) \leq 1$ for all $x \in \Re$.
- The nonlinear PDE (2.32) is not hyperbolic and is not parabolic. When $\rho(t,x) \leq 1$, we get $\rho_t(t,x) = 0$ (a zero-speed hyperbolic PDE) with the propagation speed being zero; and when $\rho(t,x) > 1$, we get the following nonlinear heat equation

$$\rho_t + h(-\kappa(\rho)\rho_x)\rho_x - \rho h'(-\kappa(\rho)\rho_x)\kappa'(\rho)\rho_x^2 = \rho h'(-\kappa(\rho)\rho_x)\kappa(\rho)\rho_{xx}. \tag{2.34}$$

- The PDE (2.32) (see also (2.34)) is a nonlinear heat equation with the diffusion coefficient depending on both the density $\rho$ and the spatial derivative of the density $\rho_x$. This feature is rarely studied in the literature, where the diffusion coefficient exclusively depends on the density $\rho$, see for instance [11], [25], [26], [27], [29], [30], [33] and references therein. Notable exceptions are [11], [20], [21], and [15], where the viscosity also depends on the spatial derivative of the state.

Define the set

$$X = \left\{ \rho \in L^1(\Re) \cap L^\infty(\Re) : \inf_{x \in \Re}(\rho(x)) = 0, \sup_{x \in \Re}(\rho(x)) < R, \lim_{x \to \pm\infty}(\rho(x)) = 0 \right\}. \tag{2.35}$$

Using (2.32), it can be shown that for every classical solution $\rho \in C^1(\Re_+ \times \Re; [0,R)) \cap C^2((0,+\infty) \times \Re)$ of (2.32), (2.33) with $\rho[t] \in X$ for all $t \geq 0$, it holds that

$$\dot{m}(t) = 0 \text{ for } t > 0 \tag{2.36}$$

where $m(t)$ is the total mass

$$m(t) = \int_{-\infty}^{+\infty} \rho(t,x) dx. \tag{2.37}$$

Therefore, the total mass remains constant. The fact that mass is conserved is of great importance, since it characterizes physically admissible solutions.

For any classical solution $\rho \in C^1(\Re_+ \times \Re; [0,R)) \cap C^2((0,+\infty) \times \Re)$ of (2.32), (2.33) with $\rho[t] \in X$

for all $t \geq 0$, we can define the following functionals

$$E_1(t) = \int_{-\infty}^{+\infty} \rho(t,x) H\left(-\kappa(\rho(t,x))\rho_x(t,x)\right) dx \tag{2.38}$$

$$E_2(t) = \int_{-\infty}^{+\infty} Q(\rho(t,x)) dx \tag{2.39}$$

where

$$H(w) := \int_0^w h(s) ds \text{ for all } w \in (-1, b) \tag{2.40}$$

$$Q(\rho) := \int_1^\rho (\rho - \tau) \kappa(\tau) d\tau \text{ for all } \rho \in [0, R]. \tag{2.41}$$

The above functionals are inspired by the mechanical energy of the original model (2.22), (2.23), with the functional $E_1$ in (2.38) expressing the kinetic energy and the functional $E_2$ in (2.39) expressing the potential energy. A direct consequence of (2.32) and the functionals defined by (2.38) and (2.39) are the following inequalities that hold for all classical solutions $\rho \in C^1(\Re_+ \times \Re; [0, R]) \cap C^2((0, +\infty) \times \Re)$ of (2.32), (2.33) with $\rho[t] \in X$ for all $t \geq 0$:

$$\dot{E}_1(t) = -\int_{-\infty}^{+\infty} \rho^2(t,x) \kappa(\rho(t,x)) \left(h(-\kappa(\rho(t,x))\rho_x(t,x))\right)_x^2 dx \leq 0, \text{ for } t > 0 \tag{2.42}$$

$$\dot{E}_2(t) = \int_{-\infty}^{+\infty} \rho(t,x) \kappa(\rho(t,x)) \rho_x(t,x) h(-\kappa(\rho(t,x))\rho_x(t,x)) dx \leq 0, \text{ for } t > 0. \tag{2.43}$$

Inequalities (2.42), (2.43) show us that both functionals (2.38), (2.39) are decreasing along classical solutions of (2.32). Inequality (2.43) gives an entropy-like condition (see (2.48) below) that characterize physically meaningful solutions, see [9], [31].

Let the initial condition

$$\rho[0] = \rho_0 \in X \tag{2.44}$$

be given. We give next the concept of a weak solution for the initial-value problem (2.32), (2.33), (2.44).

**Definition 1:** *A pair of maps* $\rho, w : \Re_+ \to L^1(\Re) \cap L^\infty(\Re)$ *with*

$$\inf_{x \in \Re} (\rho(t,x)) = 0, \text{ for all } t \geq 0, \quad \sup_{x \in \Re, t \geq 0} (\rho(t,x)) < R \tag{2.45}$$

$$-1 < \inf_{x \in \Re, t \geq 0} (w(t,x)) \leq \sup_{x \in \Re, t \geq 0} (w(t,x)) < b \tag{2.46}$$

$$\int_{-\infty}^{+\infty} \rho(t,x)dx = \int_{-\infty}^{+\infty} \rho_0(x)dx, \text{ for all } t \geq 0 \tag{2.47}$$

$$\int_{-\infty}^{+\infty} Q(\rho(t,x))dx \leq \int_{-\infty}^{+\infty} Q(\rho_0(x))dx, \text{ for all } t \geq 0 \tag{2.48}$$

*For which (2.33) holds, is a weak solution of the initial-value problem (2.32), (2.33), (2.44), if, for every function $\phi \in C^2(\Re_+ \times \Re)$ of compact support, the following equations are satisfied:*

$$-\int_{-\infty}^{+\infty} \phi(0,x)\rho_0(x)dx = \int_0^{+\infty}\int_{-\infty}^{+\infty} \rho(t,x)\big(w(t,x)\phi_x(t,x) + \phi_t(t,x)\big)dxdt \tag{2.49}$$

$$\int_0^{+\infty}\int_{-\infty}^{+\infty} \big(\phi(t,x)\beta(w(t,x)) - \phi_x(t,x)Q'(\rho(t,x))\big)dxdt = 0. \tag{2.50}$$

Definition 1 is more demanding than the typical definition of weak solution (see for instance [24], [31]) since, in addition to (2.49) and (2.50), it also requires that the mass is conserved (see (2.47)) and the entropy-like condition (2.48) also holds.

**Proposition 1:** *Suppose that $\rho \in C^1(\Re_+ \times \Re;[0,R)) \cap C^2((0,+\infty) \times \Re)$ is a classical solution of the initial-value problem (2.32), (2.33), (2.44), that satisfies (2.45) and*

$$\kappa(\rho)\rho_x \in L^\infty(\Re_+ \times \Re) \tag{2.51}$$

$$\rho[t] \in X, \text{ for all } t \geq 0. \tag{2.52}$$

*Then, the pair $\rho, w$ is also a weak solution of (2.32), (2.33), (2.44), where $w$ is given by*

$$w = h(-\kappa(\rho)\rho_x). \tag{2.53}$$

*Moreover, if the pair $\rho, w$ is a weak solution of the initial-value problem (2.32), (2.33), (2.44) with $\rho \in C^1(\Re_+ \times \Re;[0,R)) \cap C^2((0,+\infty) \times \Re)$, then $\rho$ is also a classical solution of the initial-value problem (2.32), (2.33), (2.44) that satisfies (2.51), (2.45), and (2.52).*

## 3. A Numerical Scheme and its Properties

Using finite differences for (2.32) with time-step $\delta t > 0$ and spatial discretization $\delta x > 0$, we next study the following explicit numerical scheme

$$\rho_i^+ = \rho_i + \frac{\delta t}{\delta x}(G_{i-1} - G_i), \; i \in \mathbb{Z} \tag{3.1}$$

where

$$G_i = \rho_i w_i = \rho_i h(-q_i), \quad i \in \mathbb{Z} \tag{3.2}$$

$$q_i = -\beta(w_i) = \frac{Q'(\rho_{i+1}) - Q'(\rho_i)}{\delta x}, \quad i \in \mathbb{Z} \tag{3.3}$$

with $\rho_i \in [0, R)$ for $i \in \mathbb{Z}$ and $M = \sup_{i \in \mathbb{Z}}(\rho_i) < R$ be given and $\beta: (-1, b) \to \mathfrak{R}$ being the inverse of $h: \mathfrak{R} \to (-1, b)$. Here $\rho_i$ is the numerical value of density at the point $x = i\delta x$ and at time $t \geq 0$, while $\rho_i^+$ is the numerical value of density at the point $x = i\delta x$ and at time $t + \delta t$. The numerical scheme (3.1)-(3.3) is a first-order accurate discretization and is appropriately designed to respect the entropy-like condition (2.43), for sufficiently small time-step $\delta t$, and to satisfy certain very important properties, such as conservation of mass, that are discussed in detail below.

A direct consequence of the discretization (3.1)-(3.3) above, is that the mass $m(t)$ remains constant showing that the numerical scheme is conservative, see [24]. Indeed, we have

**Proposition 2:** *For every sequence* $\rho_i \in [0, R), i \in \mathbb{Z}$ *with* $\sum_{i \in \mathbb{Z}} \rho_i < +\infty$, *it holds that* $\sum_{i \in \mathbb{Z}} \rho_i^+ = \sum_{i \in \mathbb{Z}} \rho_i$.

Notice that Proposition 2 implies that $\lim_{i \to \pm\infty}(\rho_i^+) = \inf_{i \in \mathbb{Z}}(\rho_i^+) = 0$ for every sequence $\rho_i \in [0, R)$, $i \in \mathbb{Z}$ with $\sum_{i \in \mathbb{Z}} \rho_i < +\infty$. The proof of Proposition 2 is trivial and is omitted. To complete the discretization of model (2.32), the above scheme has to be supplemented with a suitable time-step for certain properties to be preserved. In the following propositions, we obtain an upper bound on the time discretization step $\delta t$, for a given space discretization $\delta x$, which guarantees that the functional $E_2$ defined by (2.39) is non-increasing (in a discretized sense) and that (2.45) holds (again in discretized sense). In particular, it is proved that, in addition to the numerical scheme being conservative, it also respects the entropy condition (2.43).

**Proposition 3:** *Define* $f(\delta t) = \delta x \sum_{i \in \mathbb{Z}} Q(\rho_i^+)$ *for all* $\delta t \geq 0$. *Then* $\frac{df}{d\delta t}(0) = \delta x \sum_{i \in \mathbb{Z}} \rho_i q_i h(-q_i) \leq 0$.

**Proposition 4:** *Suppose that* $1 \geq \frac{\delta t}{\delta x}\left(b + 2\frac{M\|h'\|_\infty}{\delta x} \max_{0 \leq \rho \leq M}(\kappa(\rho))\right)$. *Then* $\rho_i^+ \in [0, M]$ *for all* $i \in \mathbb{Z}$, *where* $M = \sup_{i \in \mathbb{Z}}(\rho_i) < R$.

For given $\delta x$, the bound on the time-step provided in Proposition 4, establishes that $\rho_i$ is non-negative and remains bounded by $M = \sup_{i \in \mathbb{Z}}(\rho_i) < R$ for all $i \in \mathbb{Z}$. It should be noted that the condition on the time-step is an analogous nonlinear version of the well-known *Courant-Friedrichs-Lewy (CFL)* condition, see [32]. Finally, under the same bounds on time-step $\delta t$, as in Proposition 4, it is shown next that the functional $E_2$ in (2.39) is non-increasing.

**Proposition 5:** *Suppose that* $1 \geq \dfrac{\delta t}{\delta x}\left(b + 2\dfrac{M\|h'\|_\infty}{\delta x}\max_{0\leq \rho \leq M}(\kappa(\rho))\right)$. *Then*

$$\sum_{i\in\mathbb{Z}}Q(\rho_i^+) \leq \sum_{i\in\mathbb{Z}}Q(\rho_i) + \delta t\left(1 - \delta t\dfrac{4M\|h'\|_\infty}{(\delta x)^2}\max_{0\leq \rho \leq M}(\kappa(\rho))\right)\sum_{i\in\mathbb{Z}}\rho_i q_i h(-q_i)$$

$$= \sum_{i\in\mathbb{Z}}Q(\rho_i) - \delta t\left(1 - \delta t\dfrac{4M\|h'\|_\infty}{(\delta x)^2}\max_{0\leq \rho \leq M}(\kappa(\rho))\right)\sum_{i\in\mathbb{Z}}\rho_i w_i \beta(w_i) \quad (3.4)$$

$$\leq \sum_{i\in\mathbb{Z}}Q(\rho_i)$$

The following proposition provides the link between the numerical solution produced by the numerical scheme (3.1)-(3.3) and the corresponding weak solution of (2.32), (2.33), (2.44).

**Proposition 6:** *Let* $\rho_0 \in X$ *be given. Let* $\delta x, \delta t > 0$ *be given constants. Consider the sequence* $\{\rho_i^k : i, k \in \mathbb{Z}, k \geq 0\}$ *defined by the recursive formula*

$$\rho_i^{k+1} = \rho_i^k + \dfrac{\delta t}{\delta x}\left(\rho_{i-1}^k h\left(\dfrac{Q'(\rho_{i-1}^k) - Q'(\rho_i^k)}{\delta x}\right) - \rho_i^k h\left(\dfrac{Q'(\rho_i^k) - Q'(\rho_{i+1}^k)}{\delta x}\right)\right) \quad (3.5)$$

*for all* $i, k \in \mathbb{Z},\ k \geq 0$

*with*

$$\rho_i^0 = \dfrac{1}{\delta x}\int_{i\delta x}^{(i+1)\delta x}\rho_0(s)ds,\ \text{for}\ i \in \mathbb{Z}. \quad (3.6)$$

*Define the functions* $\rho_{\delta x,\delta t}, w_{\delta x,\delta t} : \Re_+ \times \Re \to \Re$ *by means of the formulas:*

$$\rho_{\delta x,\delta t}(t,x) = \rho_i^k$$

$$w_{\delta x,\delta t}(t,x) = w_i^k = h\left(\dfrac{Q'(\rho_i^k) - Q'(\rho_{i+1}^k)}{\delta x}\right) \quad (3.7)$$

*for* $x \in [i\delta x,(i+1)\delta x),\ t \in [k\delta t,(k+1)\delta t),\ i,k \in \mathbb{Z},\ k \geq 0$.

*Then for every function* $\phi \in C^2(\Re_+ \times \Re)$ *of compact support, the following equations hold:*

$$\int_{-\infty}^{+\infty}\rho_{\delta x,\delta t}(t,x)dx = \int_{-\infty}^{+\infty}\rho_0(x)dx,\ \text{for all}\ t \geq 0 \quad (3.8)$$

$$-\int_{-\infty}^{+\infty} \phi(0,x)\rho_0(x)dx = \int_0^{+\infty}\int_{-\infty}^{+\infty} \rho_{\delta x,\delta t}(t,x)\left(w_{\delta x,\delta t}(t,x)\phi_x(t,x) + \phi_t(t,x)\right)dxdt + P_1 \tag{3.9}$$

$$\int_0^{+\infty}\int_{-\infty}^{+\infty}\left(\phi(t,x)\beta\left(w_{\delta x,\delta t}(t,x)\right) - \phi_x(t,x)Q'\left(\rho_{\delta x,\delta t}(t,x)\right)\right)dxdt = P_2 \tag{3.10}$$

*where*

$$\begin{aligned}
P_1 &:= -\sum_{i\in\mathbb{Z}}\int_0^{\delta t}\int_{i\delta x}^{(i+1)\delta x}\phi_t(0,i\delta x)\rho_{\delta x,\delta t}(t,x)dxdt \\
&+\delta t\sum_{k\geq 0}\sum_{i\in\mathbb{Z}}\rho_{i-1}^k w_{i-1}^k \int_{(i-1)\delta x}^{i\delta x}\int_{(i-1)\delta x}^{\xi}\phi_{xx}(k\delta t,s)dsd\xi \\
&+\delta x\sum_{k\geq 0}\sum_{i\in\mathbb{Z}}\rho_i^{k+1}\int_{k\delta t}^{(k+1)\delta t}\int_{k\delta t}^{\tau}\phi_{tt}(l,i\delta x)dld\tau - \sum_{i\in\mathbb{Z}}\int_{i\delta x}^{(i+1)\delta x}\int_{i\delta x}^{x}\rho_0(x)\phi_x(0,l)dldx \\
&-\sum_{k\geq 0}\sum_{i\in\mathbb{Z}}\int_{k\delta t}^{(k+1)\delta t}\int_{i\delta x}^{(i+1)\delta x}\left(\int_{i\delta x}^{x}\phi_{tx}(t,s)ds + \int_{k\delta t}^{t}\phi_{tt}(\tau,i\delta x)d\tau\right)\rho_{\delta x,\delta t}(t,x)dxdt \\
&-\sum_{k\geq 0}\sum_{i\in\mathbb{Z}}\int_{k\delta t}^{(k+1)\delta t}\int_{(i-1)\delta x}^{i\delta x}\left(\int_{(i-1)\delta x}^{x}\phi_{xx}(t,s)ds + \int_{k\delta t}^{t}\phi_{xt}(\tau,(i-1)\delta x)d\tau\right)\rho_{\delta x,\delta t}(t,x)w_{\delta x,\delta t}(t,x)dxdt \\
&-\sum_{k\geq 0}\sum_{i\in\mathbb{Z}}\int_{(k+1)\delta t}^{(k+2)\delta t}\int_{i\delta x}^{(i+1)\delta x}\left(\int_{k\delta t}^{(k+1)\delta t}\phi_{tt}(\tau,i\delta x)d\tau\right)\rho_{\delta x,\delta t}(t,x)dxdt
\end{aligned} \tag{3.11}$$

*and*

$$\begin{aligned}
P_2 &:= -\sum_{k\geq 0}\sum_{i\in\mathbb{Z}}\int_{k\delta t}^{(k+1)\delta t}\int_{i\delta x}^{(i+1)\delta x}\left(\int_{i\delta x}^{x}\phi_{xx}(t,s)ds\right)Q'\left(\rho_{\delta x,\delta t}(t,x)\right)dxdt \\
&+\sum_{k\geq 0}\sum_{i\in\mathbb{Z}}\int_{k\delta t}^{(k+1)\delta t}\int_{i\delta x}^{(i+1)\delta x}\left(\int_{k\delta t}^{t}\phi_t(\tau,i\delta x)d\tau + \int_{i\delta x}^{x}\phi_x(t,\xi)d\xi\right)\beta\left(w_{\delta x,\delta t}(t,x)\right)dxdt \\
&-\sum_{k\geq 0}\sum_{i\in\mathbb{Z}}\int_{k\delta t}^{(k+1)\delta t}\int_{(i+1)\delta x}^{(i+2)\delta x}\left(\int_{i\delta x}^{(i+1)\delta x}\phi_{xx}(t,s)ds + \int_{k\delta t}^{t}\phi_{xt}(\tau,i\delta x)d\tau\right)Q'\left(\rho_{\delta x,\delta t}(t,x)\right)dxdt \\
&+\delta t\sum_{k\geq 0}\sum_{i\in\mathbb{Z}}Q'\left(\rho_{i+1}^k\right)\int_{i\delta x}^{(i+1)\delta x}\int_{i\delta x}^{\xi}\phi_{xx}(k\delta t,s)dsd\xi
\end{aligned} \tag{3.12}$$

*Furthermore, if* $1\geq \dfrac{\delta t}{\delta x}\left(b + 2\dfrac{\|\rho_0\|_\infty \|h'\|_\infty}{\delta x}\max_{0\leq \rho\leq \|\rho_0\|_\infty}(\kappa(\rho))\right)$ *and* $1\geq \delta t\dfrac{4\|\rho_0\|_\infty \|h'\|_\infty}{(\delta x)^2}\max_{0\leq \rho\leq \|\rho_0\|_\infty}(\kappa(\rho))$
*then we also have for all* $t\geq 0$:

$$\int_{-\infty}^{+\infty} Q(\rho_{\delta x,\delta t}(t,x))dx \leq \int_{-\infty}^{+\infty} Q(\rho_0(x))dx \tag{3.13}$$

$$0 = \inf_{x \in \mathbb{R}}(\rho_{\delta x,\delta t}(t,x)) \leq \sup_{x \in \mathbb{R}}(\rho_{\delta x,\delta t}(t,x)) \leq \|\rho_0\|_\infty \tag{3.14}$$

Proposition 6 is not adequate to guarantee existence of weak solutions in the sense of Definition 1. To be able to show existence of a weak solution using Proposition 6, we need to show additional properties. For instance, we may need to show the existence of two sequences: a sequence $\{(\delta x)_q > 0 : q \in \mathbb{N}\}$ which converges to zero and a sequence $\{(\delta t)_q > 0 : q \in \mathbb{N}\}$ satisfying $1 \geq \frac{(\delta t)_q}{(\delta x)_q}\left(b + 2\frac{\|\rho_0\|_\infty \|h'\|_\infty}{(\delta x)_q} \max_{0 \leq \rho \leq \|\rho_0\|_\infty}(\kappa(\rho))\right)$ and $1 \geq (\delta t)_q \frac{4\|\rho_0\|_\infty \|h'\|_\infty}{(\delta x)_q^2} \max_{0 \leq \rho \leq \|\rho_0\|_\infty}(\kappa(\rho))$ for all $q \in \mathbb{N}$ for which there exists a constant $L > 0$ (independent from $(\delta t)_q$ and $(\delta x)_q$) such that the corresponding numerical solution with $\delta t = (\delta t)_q$ and $\delta x = (\delta x)_q$ satisfies

$$\left|Q'(\rho_{i+1}^k) - Q'(\rho_i^k)\right| \leq L\delta x \text{ , for } i,k \in \mathbb{Z}, \ k \geq 0, \tag{3.15}$$

Finally, the following proposition, shows the asymptotic behavior of the numerical scheme. It also shows that the possibility of an inequality like (3.15) to be satisfied by the computed numerical solution is large.

**Proposition 7:** *Suppose that* $1 \geq \frac{\delta t}{\delta x}\left(b + 2\frac{\|\rho_0\|_\infty \|h'\|_\infty}{\delta x} \max_{0 \leq \rho \leq \|\rho_0\|_\infty}(\kappa(\rho))\right)$ *and* $1 > \delta t \frac{4\|\rho_0\|_\infty \|h'\|_\infty}{(\delta x)^2} \max_{0 \leq \rho \leq \|\rho_0\|_\infty}(\kappa(\rho))$. *Then*

$$\lim_{k \to +\infty}\left(\rho_i^k \left|Q'(\rho_i^k) - Q'(\rho_{i+1}^k)\right|\right) = 0 \text{ , for all } i \in \mathbb{Z} \tag{3.16}$$

## 4. Numerical Experiments

*4.1 An Academic Example*

The aim of this example is to demonstrate the numerical solution of model (2.32) and the convergence of the functionals defined by (2.35) and (2.38). In order to apply the numerical scheme (3.1)-(3.3), we select $x \in [-1,3]$ and the following functions:

$$h(s) = tanh(s), \quad s \in \mathfrak{R} \tag{4.1}$$

$$\kappa(\rho) = \begin{cases} 0 & ,0 \leq \rho \leq 1 \\ c\dfrac{(\rho-1)^2}{R-\rho} & ,1 < \rho < R \end{cases}. \tag{4.2}$$

The constant $c > 0$ is a parameter that allows to increase or decrease the viscosity. Notice also that: (i) $h(0) = 0$, $h' \in L^\infty(\mathfrak{R})$ and $h'(s) > 0$ for all $s \in \mathfrak{R}$; and (ii) the function $\kappa(\rho)$, defined by (4.2), is a $C^1$ function with $\kappa(\rho) = 0$ for all $\rho \in [0,1]$, $\kappa(\rho) > 0$ for all $\rho \in (1,R)$ and $\lim_{\rho \to R^-} (\kappa(\rho)) = +\infty$.

Using (4.2) and definition (2.40) we also obtain that

$$Q(\rho) = c \begin{cases} 0 & ,0 \leq \rho \leq 1 \\ (\rho-1)\left(\dfrac{\rho^2+\rho+1}{3} + (R-\rho)\dfrac{\rho+2R+1}{2} + \rho - 2R\right) \\ \quad + (R-1)^2(\rho-R)\ln\left(\dfrac{R-1}{R-\rho}\right) & ,1 < \rho < R \end{cases} \tag{4.3}$$

Let $M = 1.4$, where $M = \sup_{i \in \mathbb{Z}}(\rho_i) < R$ (recall Proposition 4), spatial-discretization $\delta x = 0.04$ and consider the initial condition:

$$\rho[0] = \begin{cases} 0.25(x-\varepsilon_1)^2(x-\varepsilon_2)^2 & ,\varepsilon_1 < x < \varepsilon_2 \\ 0 & \text{otherwise} \end{cases}$$

where $\varepsilon_1 = -0.52$ and $\varepsilon_2 = 2.52$. We also let $R = 2$, $b = 1$, and due to definition (4.1) and the fact that $\kappa(\rho)$ in (4.2) is an increasing function, we get that $\max_{0 \leq \rho \leq M}(\kappa(\rho)) = \kappa(M)$ and $\|h'\|_\infty = 1$. We also let $\delta t = 10^{-4}$, which satisfies

$$\delta t \leq \dfrac{\delta x^2}{\delta x b + 2M \|h'\|_\infty \kappa(M)} \tag{4.4}$$

for all values of $c \in \{1,5,10,15\}$. Figure 1 displays the density profile for different values of constant $c$. By increasing the value of the constant $c$, the convergence rate of the numerical solution to an equilibrium $\rho(x) \leq 1$, $w(x) \equiv 0$ is faster. This is also illustrated in Figure 2 that shows the convergence of the logarithm of the functionals $E_1(t)/E_1(0)$ and $E_2(t)/E_1(0)$.

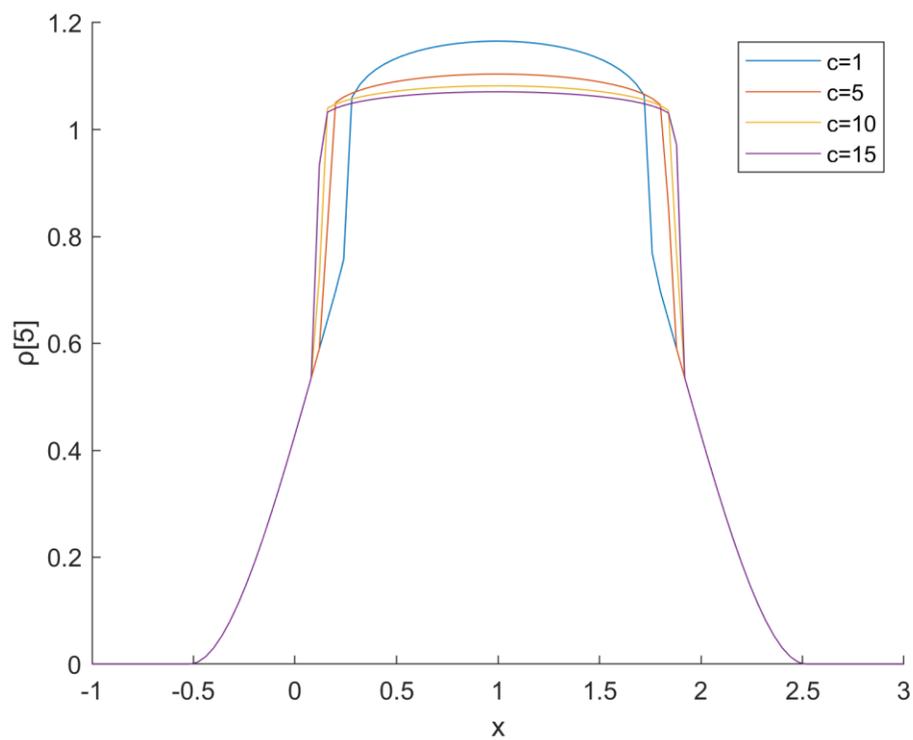

**Figure 1**: The density profile $\rho[5]$ for different values of the constant $c$.

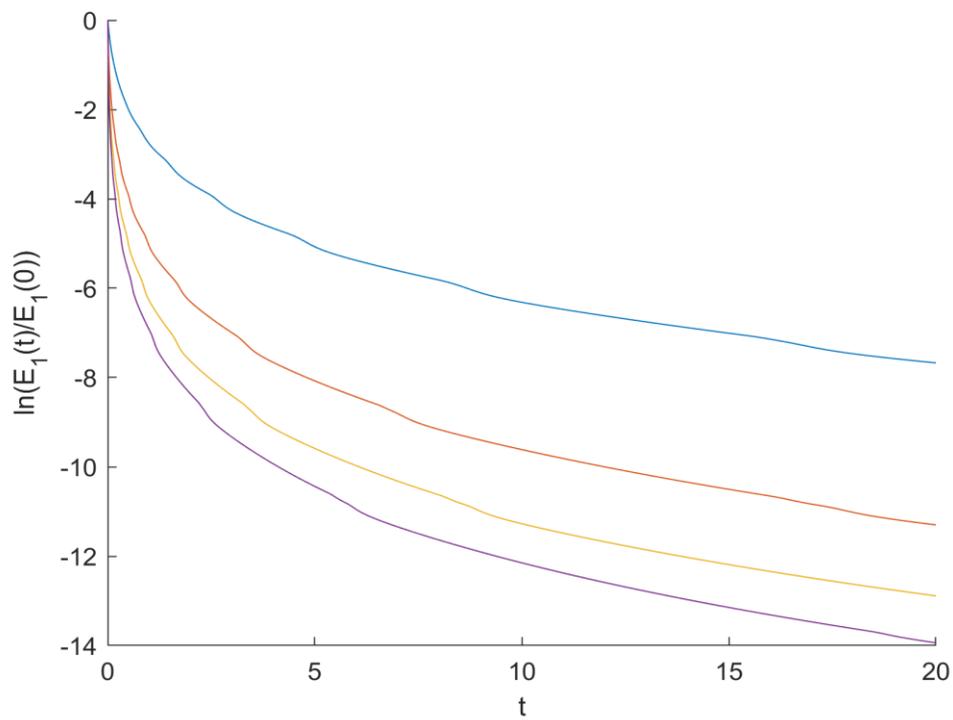

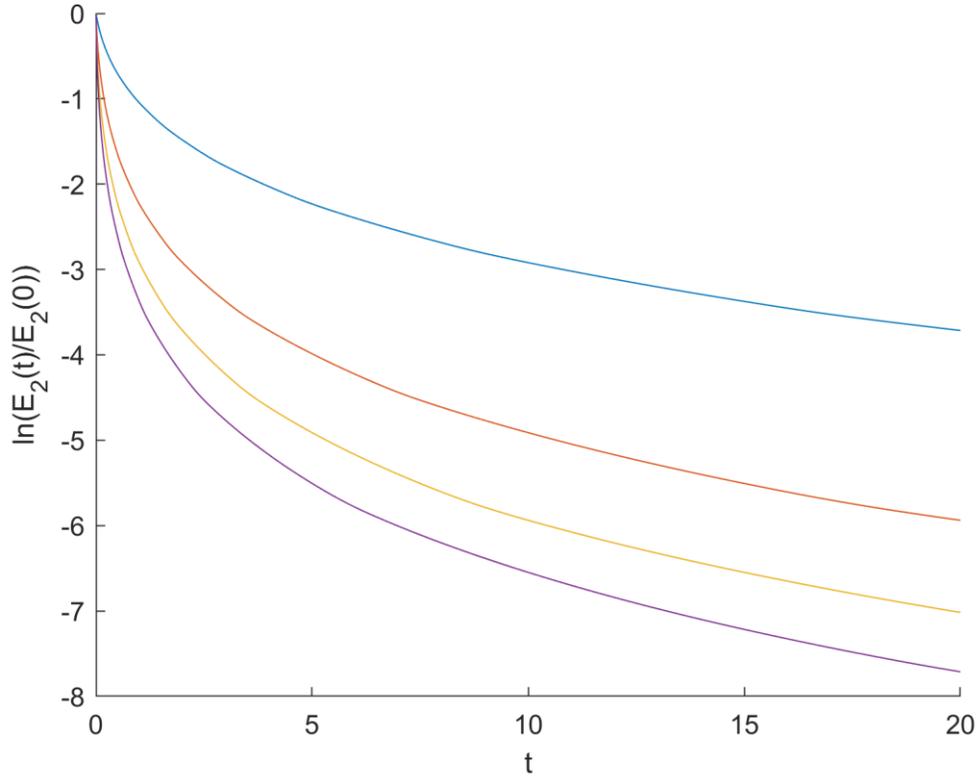

**Figure 2**: Time-evolution of $\ln(E_1(t)/E_1(0))$ and $\ln(E_2(t)/E_2(0))$
for different values of constant $c$.

In addition, Figure 3 shows the profiles $\rho[t]$ and $w[t]$ for $c=1$ and $\delta x = 0.04$, $\delta t = 10^{-3}$ at different times indicating the convergence to an equilibrium $\rho \le 1$, $w = 0$, where $w(t,x) = h(-\kappa(\rho(t,x)\rho_x(t,x)))$.

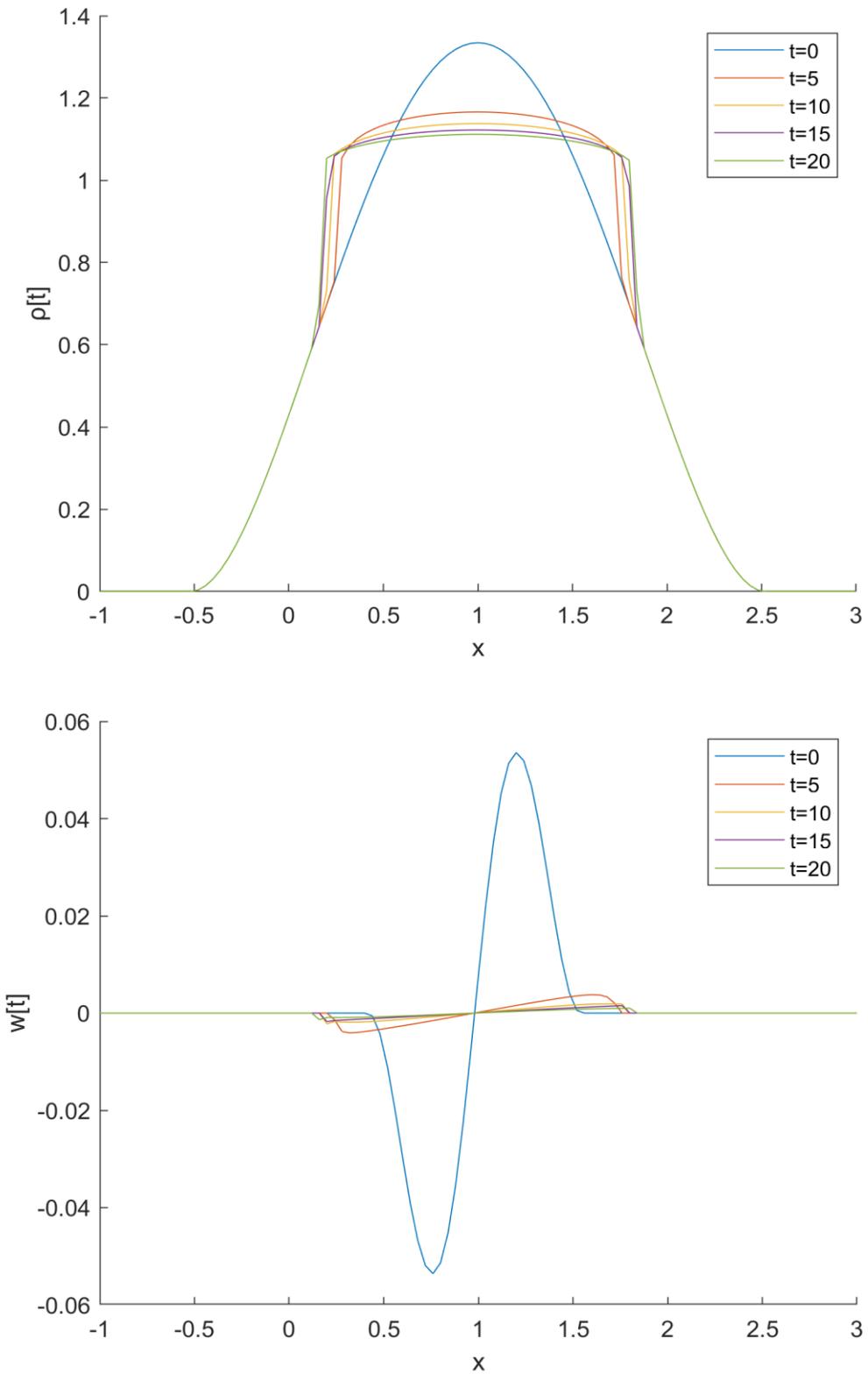

**Figure 3**: The profiles $\rho[t]$ and $w[t]$ for $c = 1$.

The numerical scheme (3.1)-(3.3) is explicit and is supplemented with a suitable time-step for the preservation of certain properties (recall Propositions 4 and 5). This restriction in the size of the time-step can be removed when implicit finite difference schemes are used. To this purpose, it is interesting to compare the solution of the explicit numerical scheme (3.1)-(3.3) with the solution of an implicit numerical scheme. The implicit analogue of numerical scheme (3.1)-(3.3) is the following:

$$\rho_i^+ = \rho_i + \frac{\delta t}{\delta x}\left(G_{i-1}^+ - G_i^+\right), \quad i \in \mathbb{Z} \tag{4.5}$$

where

$$G_i^+ = \rho_i^+ w_i^+ = \rho_i^+ h\left(-q_i^+\right), \quad i \in \mathbb{Z} \tag{4.6}$$

$$q_i^+ = -\beta(w_i^+) = \frac{Q'(\rho_{i+1}^+) - Q'(\rho_i^+)}{\delta x}, \quad i \in \mathbb{Z}. \tag{4.7}$$

In Figure 4, the sup-norm of the difference between the solution of the implicit scheme (4.5)-(4.7) and the solution of the explicit scheme (3.1)-(3.3) is illustrated. By increasing the time-step of the implicit scheme, which is represented in the horizontal axis, and using time-step $\delta t = 10^{-4}$ for the explicit scheme, the sup-norm of the difference between the solution of the implicit and the explicit schemes is also increased. Furthermore, it is obvious that $\lim_{t \to \infty} \left\|\rho_{imp}[t] - \rho_{exp}[t]\right\|_\infty = 0$ and $\lim_{t \to \infty} \left\|w_{imp}[t] - w_{exp}[t]\right\|_\infty = 0$, indicating that the solutions of the two schemes coincide as they converge to the equilibrium $\rho(x) \leq 1$, $w(x) \equiv 0$. This is also illustrated in Figure 5 which shows the sup-norm of the difference between the solution of the implicit scheme (4.5)-(4.7) using time-step $\delta t = 10^{-2}$ and the solution of the explicit scheme (3.1)-(3.3) using time-step $\delta t = 10^{-4}$.

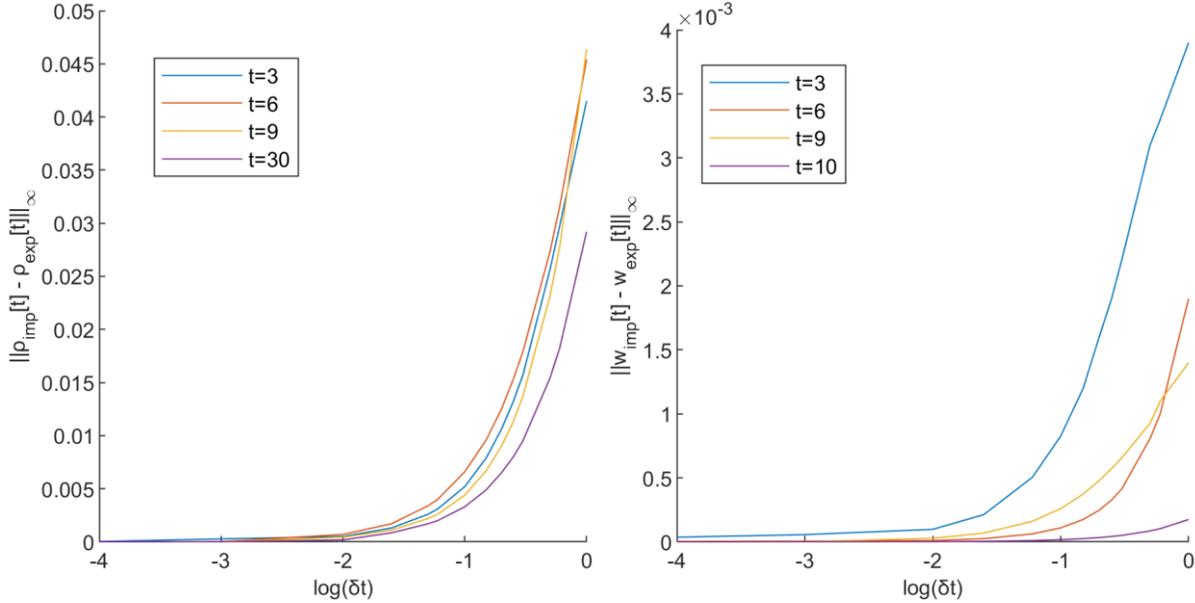

**Figure 4**: Time-evolution of $\left\|\rho_{imp}[t] - \rho_{exp}[t]\right\|_\infty$ and $\left\|w_{imp}[t] - w_{exp}[t]\right\|_\infty$ using time-step $\delta t = 10^{-4}$ for the explicit scheme and different time-steps for the implicit scheme.

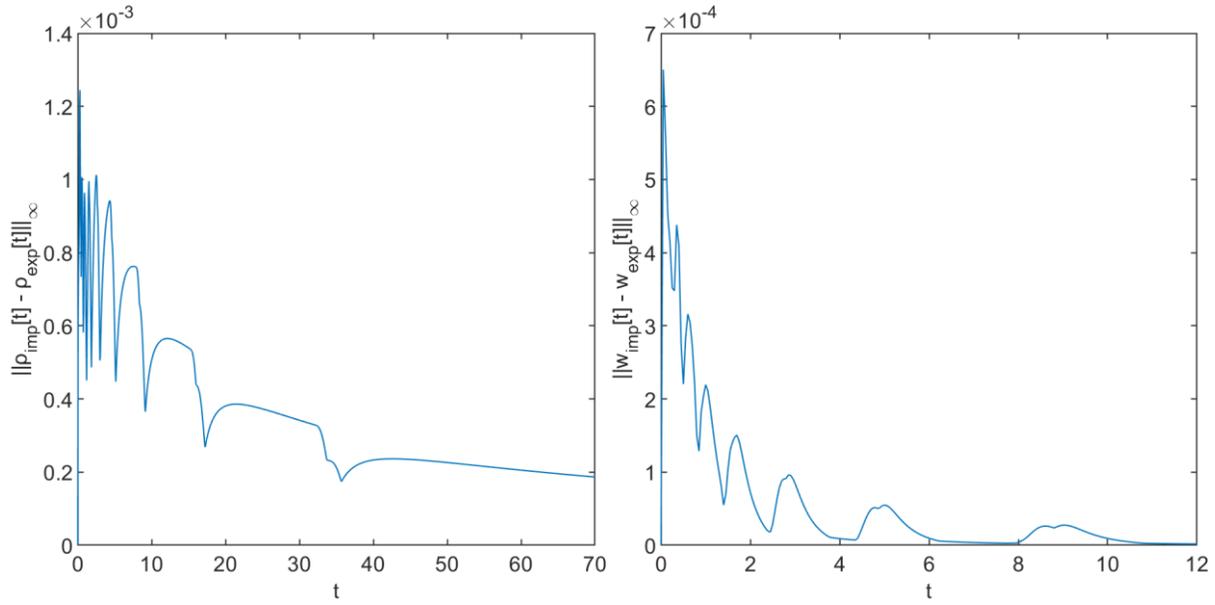

**Figure 5**: Time-evolution of $\|\rho_{imp}[t] - \rho_{exp}[t]\|_\infty$ and $\|w_{imp}[t] - w_{exp}[t]\|_\infty$ using time-step $\delta t = 10^{-2}$ for the implicit scheme and $\delta t = 10^{-4}$ for the explicit scheme.

## 4.2. A Traffic Application

In this section, we consider a traffic scenario and compare the density and flow of (2.32) with one of the most well-known traffic flow models for human drivers, the so-called LWR model ([22], [28])

$$\tilde{\rho}_\tau + (\tilde{\rho}\tilde{v})_\xi = 0 \ , \ \tau > 0, \xi \in \Re \tag{4.8}$$

with the speed $\tilde{v}$ given by

$$\tilde{v}(\tau,\xi) = v_f \exp\left(-\frac{1}{a}\left(\frac{\tilde{\rho}(\tau,\xi)}{\rho_c}\right)^a\right) \tag{4.9}$$

where $v_f > 0$ is the free-flow speed, $\rho_c$ is the critical density, and $a > 0$ is a parameter.

Our model (2.32), in the original variables $\tau > 0$, $\xi \in \Re$, is given by equation (4.8) with the speed $\tilde{v}$ defined by

$$\tilde{v}(t,\xi) = v^*\left(1 + h\left(-\mu\left(\frac{\tilde{\rho}(\tau,\xi)}{\bar{\rho}}\right)\frac{\tilde{\rho}_\xi(\tau,\xi)\bar{\rho}}{\tilde{\rho}^2(\tau,\xi)}\right)\right) \tag{4.10}$$

(recall the transformation $\xi = rx + v^*\tau$, $x \in \Re$ and the quantities $\tau = \frac{r}{v^*}t$, $w = \frac{\tilde{v} - v^*}{v^*}$, $\rho = \frac{\tilde{\rho}}{\bar{\rho}}$, with $r = 1$), where

$$h(s) = \beta^{-1}(s), \quad s \in \Re \quad (4.11)$$

with $\beta(w)$ defined by (2.19), and $\mu(\rho)$ given by

$$\mu(\rho) = \begin{cases} 0 & , 0 \leq \rho \leq 1 \\ c\frac{(\rho-1)^2}{R-\rho} & , 1 < \rho < R \end{cases} \quad (4.12)$$

where $c > 0$ is a parameter. Therefore, (4.8), (4.10), is a first-order automated vehicle model (AV model).

We consider a single-lane motorway with $v_f = 102$ $(km/h)$, $\rho_{\max} = 180$ $(veh/km)$, $\rho_c = 33.3$ $(veh/km)$, and $a = 2.34$ that were estimated in [19] based on real data from a part of the Amsterdam A10 motorway. In this scenario, we assume that the initial density on the road is characterized by a congestion belt in the interval $[1.5, 2.75]$ (km) as shown in Figure 6 with $\tilde{\rho}_0(\xi) = 0$ for $\xi < 0$ and $\xi > 0$. Note also that outside the congested area, the density is below the critical density $\rho_c$.

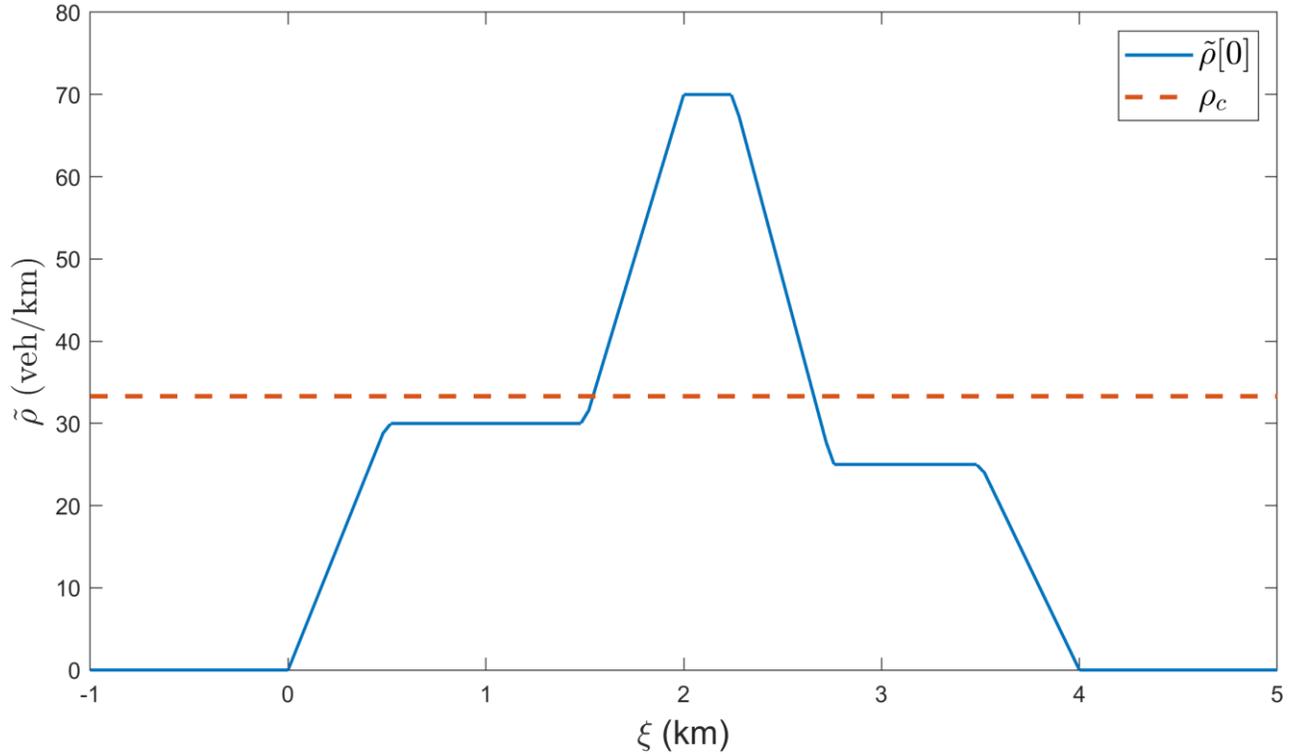

**Figure 6:** Initial density $\tilde{\rho}_0(\xi)$.

For the parameters of AV Model, we select the speed set-point $v^* = 70\ (km/h)$; significantly smaller than the speed limit $v_{max} = 110\ (km/h)$, the viscosity constant $c = 40$, interaction density $\bar{\rho} = 31$, and recall that the quantities $R$ and $b$ are given by $R = \dfrac{\rho_{max}}{\bar{\rho}} > 1$ and $b = \dfrac{v_{max} - v^*}{v^*}$, respectively. For the LWR model, we consider the interval $\xi \in [0,120]$ and use Godunov's method to obtain its solution. The density profiles for both models over a simulation time of $1h$ are shown in Figure 7.

Figure 7 shows that the density of both models dissipates along the road. However, for the LWR model the dissipation is much stronger and the density is spread over a large road interval for increasing $\tau > 0$. More specifically, for $\tau = 1\ (h)$, the density is non-zero over the interval $[90,106]$ (km), which implies that the vehicles retain large inter-vehicle distances while the maximum density is equal to $\max\limits_{\xi \in [90,106]} (\tilde{\rho}[1]) = 11.8\ (veh/km)$. On the other hand, for the AV Model, the density dissipates at a lower rate but the vehicles remain in a $4km$ stretch of the road ($\xi \in [70,74]$ (km), since the desired speed is $v^* = 70\ (km/h)$), as was the case in the initial density, and with maximum density $\max\limits_{\xi \in [70,74]} (\tilde{\rho}[1]) = 34.83\ (veh/km)$ (near the critical density $\rho_c$). Note however that the travel time of the AV model is lower since vehicles are moving with a significantly lower speed. Note also that the density for AV Model converges towards an equilibrium where $\tilde{v}(\xi) \equiv v^*$ and $\tilde{\rho}(\xi) \leq \bar{\rho} < \rho_c$, see also Figure 8 which displays the speed for both models.

Moreover, we also consider the case where the speed set-point is equal to the free-flow speed $v^* = v_f = 102\ (km/h)$. The density profiles of the AV model with $v^* = 102\ (km/h)$ are shown in Figure 8. Note that when vehicles are moving with free-flow speed $v^* = 102\ (km/h)$, the overall travel time is higher with the AV model, since, for $\tau = 1\ (h)$, vehicles are in the interval $(102,106)$ (km), whereas with the LWR, vehicles are in the interval $(90,106)$ (km), compare with Figure 7 (top).

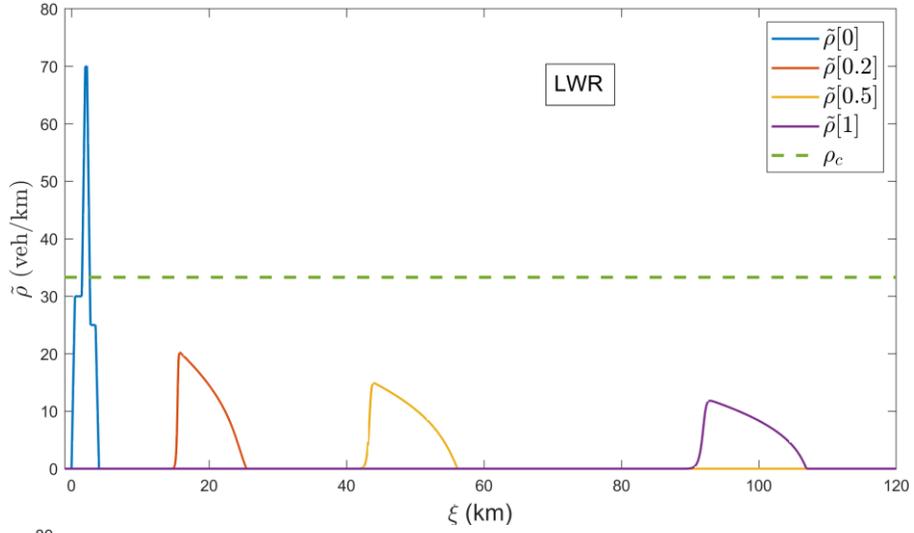

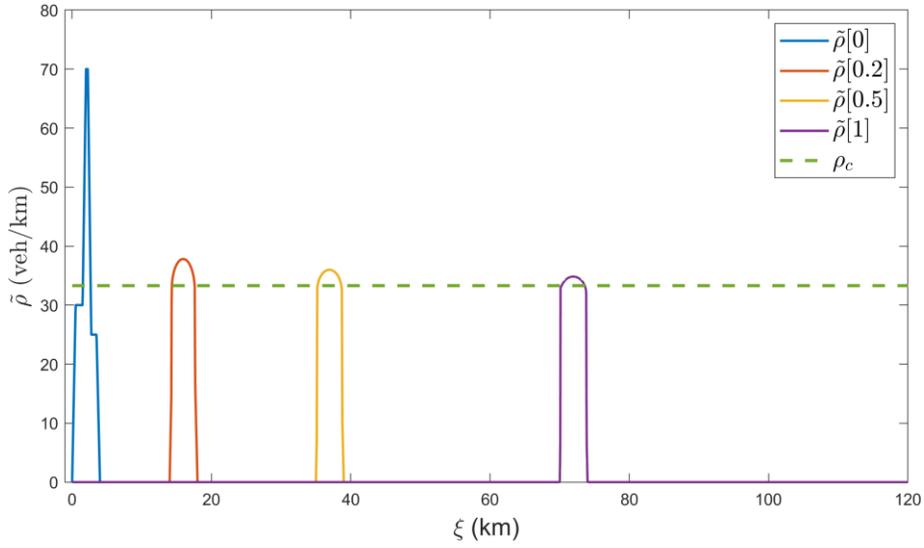

**Figure 7:** Density profiles for the LWR (top), and density profiles of AV Model (bottom).

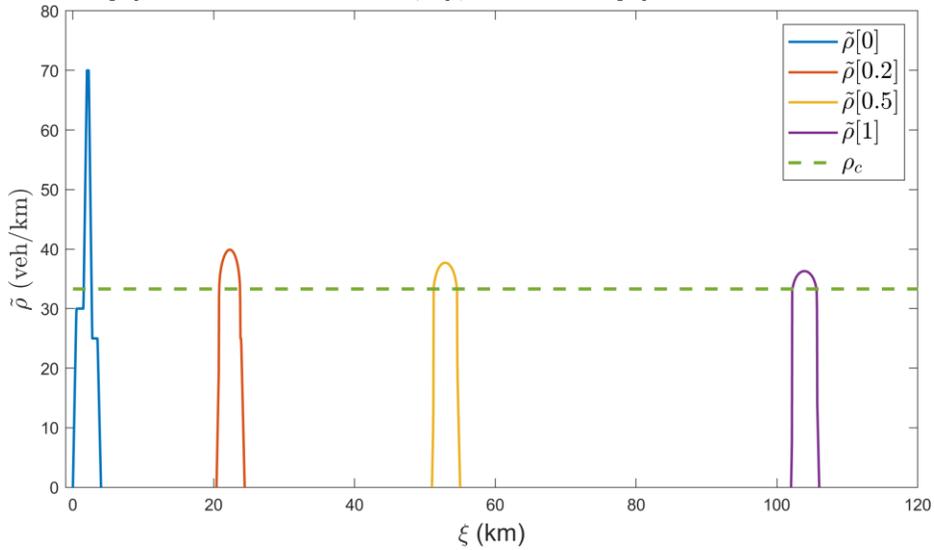

**Figure 8:** Density profiles of AV Model with $v^* = 102$ ($km/h$)

Since with the AV model (4.8), (4.10), the density converges towards the critical density while vehicles are spread in a small space interval (recall Figure 7), it is expected that the flow $q = \tilde{\rho}\tilde{v}$ of the AV model (4.8), (4.10) is much higher than the flow of the LWR model. Equivalently, vehicles equipped with the cruise controllers proposed in [16] will retain smaller inter-vehicle distances (higher density) than human-driven vehicle. Indeed, to compare the mean flow for each model, we define

$$\text{Mean Flow} = \frac{1}{T}\int_0^T \int_{a(t)}^{b(t)} \frac{\tilde{\rho}(t,x)\tilde{v}(t,x)}{b(t)-a(t)} dx dt$$

where $T > 0$ denotes the maximum time and $[a(\tau), b(\tau)]$ denotes the interval where $\tilde{\rho}(\tau, \xi) \neq 0$ for $\xi \in [a(\tau), b(\tau)]$ at each time instant $\tau \geq 0$. We consider the cases: (i) $v^* = 70$ (km/h) – the speed set point is much lower than the free flow speed; and (ii) $v^* = 102$ (km\h) – the speed set point is equal to the free-flow speed $v_f$. The mean flow for both models and $T = 1$ $(h)$ is shown in Table 1. In both cases, the mean flow for the AV model (4.8), (4.10) is much higher than the mean flow of the LWR. In particular, when the vehicles are moving with free-flow speed, the mean flow of the AV model is 346% higher than the mean flow of the LWR model. The density profiles of the AV model with $v^* = 102$ (km\h) are shown in Figure 9. Note that when vehicles are moving with free-flow speed, the overall travel time is higher with the AV model, since for $\tau = 1$ (h), vehicles are in the interval $\xi \in (102, 106)$ (km) whereas with the LWR, vehicles are in the interval $\xi \in (90, 106)$ (km).

**Table 1**

|  | AV Model (4.8), (4.10) with $v^* = 70$ (km/h) | AV Model (4.8), (4.10) with $v^* = 102$ (km/h) | LWR Model (4.8), (4.9) with $v_f = 102$ (km/h) |
|---|---|---|---|
| Mean Flow | 2232 (veh/h) | 3201 (veh/h) | 926 (veh/h) |

Finally, it should be noticed that the macroscopic model (2.1), (2.2) follows directly from the corresponding cruise controllers in [16], whose design has many degrees of freedom and can influence the characteristics of (2.1), (2.2). For instance, just by increasing (or decreasing) the value of the viscosity constant $c > 0$, we can influence the dissipation rate of the density $\tilde{\rho}$ and speed $\tilde{v}$ to an equilibrium (recall Figure 1).

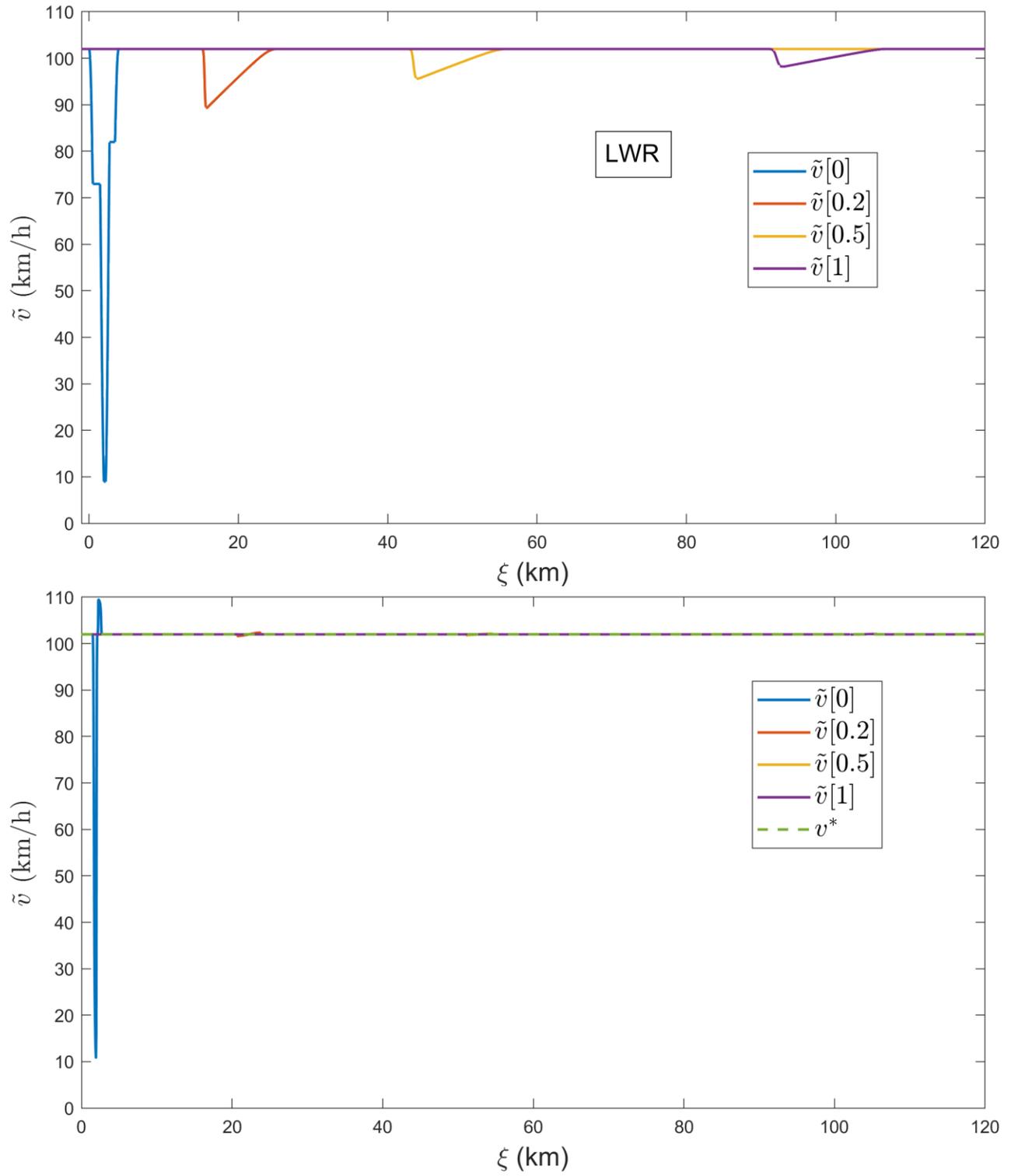

**Figure 9:** Speed profiles for the LWR (top), and the speed profiles of AV Model (bottom) with $v^* = 102$ (km/h).

# 5. Conclusions

We have presented a novel nonlinear heat equation arising from automated traffic flow models and we have shown that the solution of the nonlinear heat equation approximates the solution of the traffic model. Moreover, a concept of weak solution was defined and certain entropy-like conditions were derived. A new first-order finite difference scheme was also proposed that respects the corresponding entropy conditions and certain links between the weak solution and the numerical scheme were provided. Finally, a traffic simulation scenario and a comparison with the Lighthill-Witham-Richards (LWR) model were provided, illustrating the benefits of the proposed traffic model.

# Appendix: Proofs

**Proof of Proposition 1:** We first show that if $\rho \in C^1(\Re_+ \times \Re;[0,R]) \cap C^2((0,+\infty) \times \Re)$ is a classical solution of (2.32), (2.33), (2.44) that satisfies (2.51), (2.45), and (2.52), then the pair $\rho, w$ with $w$ defined by (2.53), is also a weak solution of the initial-value problem (2.32), (2.33), (2.44). First notice that by virtue of (2.52), (2.37), (2.36), (2.39) and (2.43), it follows that (2.47), (2.48) hold. Furthermore, it follows from (2.51), definition (2.53) and the fact that $h(\Re) = (-1, b)$, $\kappa(\rho)\rho_x \in L^\infty(\Re_+ \times \Re)$ that (2.46) holds. Equation (2.33) is a direct consequence of (2.52) and definition (2.35).

Let $\phi \in C^2(\Re_+ \times \Re)$ be a function with compact support in $\Omega := [0,T] \times [-a,a]$, for $T, a > 0$ sufficiently large. Multiplying (2.32) by $\phi$ and using (2.53), we obtain

$$\int_0^{+\infty}\int_{-\infty}^{+\infty} \phi(t,x)\left(\rho_t(t,x) + (w(t,x)\rho(t,x))_x\right)dxdt = 0 \tag{5.1}$$

Integrating by parts (5.1) and using (2.45), (2.46), Fubini's Theorem, (2.44) and that $\phi$ has compact support in $\Omega$ we obtain

$$0 = \int_0^{+\infty}\int_{-\infty}^{+\infty} \phi(t,x)\left(\rho_t(t,x) + (w(t,x)\rho(t,x))_x\right)dtdx = \int_{-\infty}^{+\infty} \rho(t,x)\phi(t,x)\Big|_{t=0}^{t=+\infty} dx - \int_{-\infty}^{+\infty}\int_0^{+\infty} \phi_t(t,x)\rho(t,x)dtdx$$
$$+ \int_0^{+\infty}\left(\phi(t,x)w(t,x)\rho(t,x)\Big|_{x=-\infty}^{x=+\infty}\right)dt - \int_0^{+\infty}\int_{-\infty}^{+\infty} \rho(t,x)w(t,x)\phi_x(t,x)dxdt$$
$$= -\int_{-\infty}^{+\infty}\phi(0,x)\rho_0(x)dx - \int_0^{+\infty}\int_{-\infty}^{+\infty}\left(\rho(t,x)w(t,x)\phi_x(t,x)\right)dxdt - \int_0^{+\infty}\int_{-\infty}^{+\infty}\phi_t(t,x)\rho(t,x)dxdt$$

(5.2)

which implies (2.49). Notice that due to definition (2.41), and the fact that $\beta$ is the inverse of $h$, (2.53) is written as follows:

$$\beta(w(t,x)) + (Q'(\rho(t,x)))_x = 0 \tag{5.3}$$

Using (5.3) and integrating by parts the equation

$$\int_0^{+\infty}\!\!\int_{-\infty}^{+\infty} \phi(t,x)\big(\beta(w(t,x))+\big(Q'(\rho(t,x))\big)_x\big)dxdt = 0 \tag{5.4}$$

we obtain

$$\int_0^{+\infty}\!\!\int_{-\infty}^{+\infty} \phi(t,x)\beta(w(t,x))-\phi_x(t,x)Q'(\rho(t,x))dxdt + \int_0^{+\infty}\phi(t,x)Q'(\rho(t,x))\big|_{x=-\infty}^{x=+\infty}dt = 0$$

and due to the fact that $\phi$ has compact support in $\Omega$, equation (2.50) follows, establishing that the pair $\rho,w$ is a weak solution.

Conversely, suppose that the pair $\rho,w$ is a weak solution of (2.32), (2.33), (2.44) as per Definition 1. Since $\rho \in C^1(\Re_+ \times \Re;[0,R)) \cap C^2((0,+\infty)\times \Re)$, it follows by definition (2.41) that $Q'(\rho) \in C^1(\Re_+ \times \Re)$.

Let arbitrary $T, a > 0$ be given. Define

$$\Omega := (0,T)\times(-a,a) \tag{5.5}$$

Let an arbitrary function $\phi \in C^2(\Re_+ \times \Re)$ with compact support in $\Omega$. From (2.50) and using integration by parts we obtain the following equation

$$\begin{aligned}
0 &= \iint_\Omega \big(\phi(t,x)\beta(w(t,x))-\phi_x(t,x)Q'(\rho(t,x))\big)dxdt \\
&= \iint_\Omega \phi(t,x)\big(\beta(w(t,x))+\big(Q'(\rho(t,x))\big)_x\big)dxdt + \int_0^T \phi(t,x)Q'(\rho(t,x))\big|_{x=-a}^{x=a}dt \\
&= \iint_\Omega \phi(t,x)\big(\beta(w(t,x))+\big(Q'(\rho(t,x))\big)_x\big)dxdt
\end{aligned} \tag{5.6}$$

Moreover, compactness of $\bar\Omega$ and (2.46) (which imply that $\beta(w)\in L^\infty(\Omega)\subseteq L^2(\Omega)$), definition (2.41) (which implies that $\big(Q'(\rho)\big)_x = \kappa(\rho)\rho_x$) and (2.51) (which implies that $\big(Q'(\rho)\big)_x \in L^\infty(\Omega)\subseteq L^1(\Omega)$) give that $\beta(w)+\big(Q'(\rho)\big)_x \in L^1(\Omega)$. Corollary 4.24 on page 110 in [5] gives that

$$\beta(w(t,x))+\big(Q'(\rho(t,x))\big)_x = 0 \text{ a.e. in } \Omega$$

or (using (2.41))

$$w(t,x) = h\big(-\kappa(\rho(t,x))\rho_x(t,x)\big) \text{ a.e. in } \Omega$$

where $h$ is the inverse of $\beta$. Finally, due to continuity of $h\big(-\kappa(\rho)\rho_x\big)$, we get that $w \in C^0(\bar\Omega)$ and using (5.5) we get

$$w(t,x) = h\left(-\kappa(\rho(t,x))\rho_x(t,x)\right) \quad \text{in } \bar{\Omega} = [0,T] \times [-a,a] \tag{5.7}$$

Since $\rho \in C^1(\mathfrak{R}_+ \times \mathfrak{R};[0,R)) \cap C^2((0,+\infty) \times \mathfrak{R})$ and since $T, a > 0$ are arbitrary, we obtain that $w \in C^0(\mathfrak{R}_+ \times \mathfrak{R};(-1,b)) \cap C^1((0,+\infty) \times \mathfrak{R})$ and (2.53) holds for all $(t,x) \in \mathbb{R}_+ \times \mathbb{R}$.

Let arbitrary $a > 0$, $T > t_0 > 0$ be given. Define

$$\Omega := (t_0, T) \times (-a, a) \tag{5.8}$$

Next, let arbitrary function $\phi \in C^2(\mathfrak{R}_+ \times \mathfrak{R})$ with compact support in $\Omega$ be given. Since $\rho \in C^1(\mathfrak{R}_+ \times \mathfrak{R};[0,R)) \cap C^2((0,+\infty) \times \mathfrak{R})$, $w \in C^0(\mathfrak{R}_+ \times \mathfrak{R};(-1,b)) \cap C^1((0,+\infty) \times \mathfrak{R})$ and since $\phi$ has compact support in $\Omega$, it follows from equations (2.49) and (2.53) by direct integration by parts that the following equation holds

$$\iint_\Omega \phi(t,x)\left(\rho_t(t,x) + (\rho(t,x)w(t,x))_x\right) dx\, dt = 0 \tag{5.9}$$

Since $\rho_t + (\rho w)_x \in C^0(\bar{\Omega}) \subset L^1(\Omega)$, Corollary 4.24 on page 110 in [5] gives that

$$\rho_t(t,x) + (\rho w)_x(t,x) = 0 \text{ a.e. in } \Omega$$

Since $\rho_t + (\rho w)_x \in C^0(\bar{\Omega})$, we get from (5.8):

$$\rho_t(t,x) + (\rho w)_x(t,x) = 0 \text{ in } \bar{\Omega} = [t_0, T] \times [-a, a]$$

Since $a > 0$, $T > t_0 > 0$ are arbitrary we obtain that $\rho_t + (\rho w)_x = 0$ for all $t > 0$ and $x \in \mathfrak{R}$. Consequently, using (2.53) we conclude that (2.32) holds for all $t > 0$ and $x \in \mathfrak{R}$.

Inclusion (2.51) is a direct consequence of (2.46) and (2.53). Inclusion (2.52) is a direct consequence of definition (2.35), (2.45), (2.33) and (2.47).

We complete the proof by showing that (2.44) holds. Using (2.49), Fubini's theorem, integration by parts and the fact that $\rho_t + (\rho w)_x = 0$ for all $t > 0$ and $x \in \mathfrak{R}$, we get for every function $\phi \in C^2(\mathfrak{R}_+ \times \mathfrak{R})$ of compact support:

$$-\int_{-\infty}^{+\infty}\phi(0,x)\rho_0(x)dx = \int_0^{+\infty}\int_{-\infty}^{+\infty}\rho(t,x)\big(w(t,x)\phi_x(t,x)+\phi_t(t,x)\big)dxdt$$

$$=\int_{-\infty}^{+\infty}\int_0^{+\infty}\rho(t,x)\phi_t(t,x)dtdx+\int_0^{+\infty}\int_{-\infty}^{+\infty}\rho(t,x)w(t,x)\phi_x(t,x)dxdt$$

$$=-\int_{-\infty}^{+\infty}\rho(0,x)\phi(0,x)dx-\int_{-\infty}^{+\infty}\int_0^{+\infty}\rho_t(t,x)\phi(t,x)dtdx-\int_0^{+\infty}\int_{-\infty}^{+\infty}\big(\rho(t,x)w(t,x)\big)_x\phi(t,x)dxdt$$

$$=-\int_{-\infty}^{+\infty}\rho(0,x)\phi(0,x)dx-\int_0^{+\infty}\int_{-\infty}^{+\infty}\big(\rho_t(t,x)+\big(\rho(t,x)w(t,x)\big)_x\big)\phi(t,x)dxdt$$

$$=-\int_{-\infty}^{+\infty}\rho(0,x)\phi(0,x)dx$$

Thus, we conclude that for every function $\phi \in C^2(\mathfrak{R})$ of compact support it holds that:

$$\int_{-\infty}^{+\infty}\phi(x)\big(\rho(0,x)-\rho_0(x)\big)dx = 0 \tag{5.10}$$

Using (5.10), Corollary 4.24 on page 110 in [5], the fact that $\rho_0 \in X \subset L^1(\mathbb{R})\cap L^\infty(\mathbb{R})$ and the fact that $\rho[0] \in X \subset L^1(\mathbb{R})\cap L^\infty(\mathbb{R})$ we get that $\rho(0,x)=\rho_0(x)$ for $x\in\mathbb{R}$ a.e.. Since $\rho[0]\in C^1(\mathbb{R})$, it follows that $\rho_0 \in C^1(\mathbb{R})$ and that $\rho(0,x)=\rho_0(x)$ for all $x\in\mathbb{R}$. Therefore, (2.44) holds. The proof is complete. ◁

**Proof of Proposition 3**: Notice first that from (3.1) we have that

$$\frac{d\rho_i^+}{d\delta t} = \frac{G_{i-1}-G_i}{\delta x}, \quad i\in\mathbb{Z} \tag{5.11}$$

Using (3.1) and (5.11) and definition of $f$ we get

$$\frac{df}{d\delta t}(\delta t) = \delta x\sum_{i\in\mathbb{Z}}Q'(\rho_i^+)\frac{d\rho_i^+}{d\delta t} = \sum_{i\in\mathbb{Z}}Q'(\rho_i^+)(G_{i-1}-G_i) = \sum_{i\in\mathbb{Z}}Q'(\rho_{i+1}^+)G_i - \sum_{i\in\mathbb{Z}}Q'(\rho_i^+)G_i$$

$$=\sum_{i\in\mathbb{Z}}\big(Q'(\rho_{i+1}^+)-Q'(\rho_i^+)\big)G_i = \delta x\sum_{i\in\mathbb{Z}}\rho_i\frac{Q'(\rho_{i+1}^+)-Q'(\rho_i^+)}{\delta x}h(-q_i) \tag{5.12}$$

The desired inequality is a direct consequence of (5.12), definitions (3.2), (3.3), and the fact that $h(\cdot)$ is increasing with $h(0)=0$. This completes the proof. ◁

**Proof of Proposition 4**: From (3.1) and definitions (3.2) and (3.3) we have for $i\in\mathbb{Z}$

$$\rho_i^+ = \rho_i + \frac{\delta t}{\delta x}\left(\rho_{i-1}h\left(\frac{Q'(\rho_{i-1}) - Q'(\rho_i)}{\delta x}\right) - \rho_i h\left(\frac{Q'(\rho_i) - Q'(\rho_{i+1})}{\delta x}\right)\right) \quad (5.13)$$

It follows from (2.41), (5.13), the facts that $M = \sup_{i \in \mathbb{Z}}(\rho_i)$ and $h' \in L^\infty(\mathfrak{R})$, and assumption $1 \geq \frac{\delta t}{\delta x}\left(b + 2\frac{M\|h'\|_\infty}{\delta x}\max_{0 \leq \rho \leq M}(\kappa(\rho))\right)$ that

$$\begin{aligned}
\frac{\partial \rho_i^+}{\partial \rho_i} &= 1 - \frac{\delta t}{\delta x}h\left(\frac{Q'(\rho_i) - Q'(\rho_{i+1})}{\delta x}\right) \\
&\quad - \frac{\delta t}{\delta x}\frac{\kappa(\rho_i)}{\delta x}\left(\rho_{i-1}h'\left(\frac{Q'(\rho_{i-1}) - Q'(\rho_i)}{\delta x}\right) + \rho_i h'\left(\frac{Q'(\rho_i) - Q'(\rho_{i+1})}{\delta x}\right)\right) \\
&\geq 1 - \frac{\delta t}{\delta x}b - \delta t\frac{\kappa(\rho_i)}{(\delta x)^2}\left(\rho_{i-1}h'\left(\frac{Q'(\rho_{i-1}) - Q'(\rho_i)}{\delta x}\right) + \rho_i h'\left(\frac{Q'(\rho_i) - Q'(\rho_{i+1})}{\delta x}\right)\right) \\
&\geq 1 - \frac{\delta t}{\delta x}b - 2\delta t\frac{M\|h'\|_\infty}{(\delta x)^2}\max_{0 \leq \rho \leq M}(\kappa(\rho)) \geq 0
\end{aligned} \quad (5.14)$$

Inequality (5.14) implies that $\rho_i^+$ is increasing with respect to $\rho_i$. Therefore, by the definition of $M$ and (5.13) we get that the following inequalities hold for $i \in \mathbb{Z}$

$$\rho_i^+ \leq M + \frac{\delta t}{\delta x}\left(\rho_{i-1}h\left(\frac{Q'(\rho_{i-1}) - Q'(M)}{\delta x}\right) - Mh\left(\frac{Q'(M) - Q'(\rho_{i+1})}{\delta x}\right)\right) \quad (5.15)$$

$$\rho_i^+ \geq \frac{\delta t}{\delta x}\rho_{i-1}h\left(\frac{Q'(\rho_{i-1})}{\delta x}\right) \quad (5.16)$$

Since $Q'$ is a non-decreasing function (recall definition (2.41)) and since $\rho_{i-1} \leq M$, it follows that $Q'(\rho_{i-1}) \leq Q'(M)$ and consequently $h\left(\frac{Q'(\rho_{i-1}) - Q'(M)}{\delta x}\right) \leq 0$. Moreover, since $\rho_{i+1} \leq M$, it follows that $Q'(\rho_{i+1}) \leq Q'(M)$ and consequently $h\left(\frac{Q'(M) - Q'(\rho_{i+1})}{\delta x}\right) \geq 0$. Thus, (5.15) implies that $\rho_i^+ \leq M$. Moreover, $Q'(\rho_{i-1}) \geq 0$ and consequently, $h\left(\frac{Q'(\rho_{i-1})}{2\delta x}\right) \geq 0$. Hence, from (5.16) we get $\rho_i^+ \geq 0$. This completes the proof. ◁

**Proof of Proposition 5:** Recall that $f(\delta t) = \delta x \sum_{i \in \mathbb{Z}} Q(\rho_i^+)$. Using definitions (2.6), (3.1), (3.2), equation (5.11), and inequality $\rho_i \leq M$, $i \in \mathbb{Z}$ we obtain

$$\frac{d^2 f}{d(\delta t)^2}(\delta t) = \sum_{i\in\mathbb{Z}} \rho_i \kappa(\rho_{i+1}^+) h(-q_i) \frac{d\rho_{i+1}^+}{d\delta t} - \sum_{i\in\mathbb{Z}} \rho_i \kappa(\rho_i^+) h(-q_i) \frac{d\rho_i^+}{d\delta t}$$

$$= \sum_{i\in\mathbb{Z}} \rho_i \kappa(\rho_{i+1}^+) h(-q_i) \frac{G_i - G_{i+1}}{\delta x} - \sum_{i\in\mathbb{Z}} \rho_i \kappa(\rho_i^+) h(-q_i) \frac{G_{i-1} - G_i}{\delta x}$$

$$= \sum_{i\in\mathbb{Z}} \kappa(\rho_i^+) \frac{G_{i-1} - G_i}{\delta x} \big( \rho_{i-1} h(-q_{i-1}) - \rho_i h(-q_i) \big)$$

$$= \frac{1}{\delta x} \sum_{i\in\mathbb{Z}} \kappa(\rho_i^+) (G_{i-1} - G_i)^2 \leq \frac{1}{\delta x} \max_{0\leq\rho\leq M} (\kappa(\rho)) \sum_{i\in\mathbb{Z}} (G_{i-1} - G_i)^2 \quad (5.17)$$

$$\leq \frac{4}{\delta x} \max_{0\leq\rho\leq M} (\kappa(\rho)) \sum_{i\in\mathbb{Z}} G_i^2 = \frac{4}{\delta x} \max_{0\leq\rho\leq M} (\kappa(\rho)) \sum_{i\in\mathbb{Z}} \rho_i^2 h^2(-q_i)$$

$$\leq \frac{4M}{\delta x} \max_{0\leq\rho\leq M} (\kappa(\rho)) \sum_{i\in\mathbb{Z}} \rho_i h^2(-q_i) \leq \frac{4M \|h'\|_\infty}{\delta x} \max_{0\leq\rho\leq M} (\kappa(\rho)) \sum_{i\in\mathbb{Z}} \rho_i |q_i| |h(-q_i)|$$

$$= \frac{4M \|h'\|_\infty}{\delta x} \max_{0\leq\rho\leq M} (\kappa(\rho)) \sum_{i\in\mathbb{Z}} \rho_i (-q_i) h(-q_i) = -\frac{4M \|h'\|_\infty}{\delta x} \max_{0\leq\rho\leq M} (\kappa(\rho)) \sum_{i\in\mathbb{Z}} \rho_i q_i h(-q_i)$$

From Proposition 3 and (5.17) we finally obtain

$$f(\delta t) \leq f(0) + \delta t \left( 1 - \delta t \frac{4M \|h'\|_\infty}{(\delta x)^2} \max_{0\leq\rho\leq M} (\kappa(\rho)) \right) \delta x \sum_{i\in\mathbb{Z}} \rho_i q_i h(-q_i)$$

which implies inequality (3.4). The proof is complete. ◁

**Proof of Proposition 6:** For any $\phi \in C^2(\mathfrak{R}_+ \times \mathfrak{R})$ with compact support, define $\phi_i^k = \phi(k\delta t, i\delta x)$ for $i,k \in \mathbb{Z}$, $k \geq 0$. We have from (3.5) and definition (3.7) after adding and subtracting terms that

$$-\sum_{i\in\mathbb{Z}} \frac{\rho_i^0 \phi_i^0}{\delta t} = \sum_{k\geq 0} \sum_{i\in\mathbb{Z}} \frac{\rho_{i-1}^k w_{i-1}^k}{\delta x} (\phi_i^k - \phi_{i-1}^k) + \sum_{k\geq 0} \sum_{i\in\mathbb{Z}} \rho_i^{k+1} \frac{\phi_i^{k+1} - \phi_i^k}{\delta t} \quad (5.18)$$

Moreover, from (5.18), (3.6), and (3.7) we obtain

$$-\sum_{i\in\mathbb{Z}} \phi_i^0 \int_{i\delta x}^{(i+1)\delta x} \rho_0(s) ds = \sum_{k\geq 0} \sum_{i\in\mathbb{Z}} \frac{\phi_i^k - \phi_{i-1}^k}{\delta x} \int_{k\delta t}^{(k+1)\delta t} \int_{(i-1)\delta x}^{i\delta x} \rho_{\delta x,\delta t}(t,x) w_{\delta x,\delta t}(t,x) dx dt$$

$$+ \sum_{k\geq 0} \sum_{i\in\mathbb{Z}} \frac{\phi_i^{k+1} - \phi_i^k}{\delta t} \int_{(k+1)\delta t}^{(k+2)\delta t} \int_{i\delta x}^{(i+1)\delta x} \rho_{\delta x,\delta t}(t,x) dx dt \quad (5.19)$$

Notice that since $\phi \in C^2(\mathfrak{R}_+ \times \mathfrak{R})$ we have that the following equations hold for $k,i \in \mathbb{Z}$, $k > 0$

$$\frac{\phi_i^{k+1} - \phi_i^k}{\delta t} = \phi_t(k\delta t, i\delta x) + \frac{1}{\delta t} \int_{k\delta t}^{(k+1)\delta t} \int_{k\delta t}^{\tau} \phi_{tt}(l, i\delta x) \, dl \, d\tau$$

$$\frac{\phi_i^k - \phi_{i-1}^k}{\delta x} = \phi_x(k\delta t, (i-1)\delta x) + \frac{1}{\delta x} \int_{(i-1)\delta x}^{i\delta x} \int_{(i-1)\delta x}^{\xi} \phi_{xx}(k\delta t, s) \, ds \, d\xi \qquad (5.20)$$

Combining (5.20) with (5.19) we get

$$\begin{aligned}
&-\sum_{i\in\mathbb{Z}} \int_{i\delta x}^{(i+1)\delta x} \phi(0,x)\rho_0(x)\,dx = \sum_{k\geq 0}\sum_{i\in\mathbb{Z}} \int_{k\delta t}^{(k+1)\delta t} \int_{(i-1)\delta x}^{i\delta x} \phi_x(t,(i-1)\delta x)\rho_{\delta x,\delta t}(t,x)w_{\delta x,\delta t}(t,x)\,dx\,dt \\
&+\sum_{k\geq 0}\sum_{i\in\mathbb{Z}} \int_{k\delta t}^{(k+1)\delta t} \int_{i\delta x}^{(i+1)\delta x} \phi_t(k\delta t,i\delta x)\rho_{\delta x,\delta t}(t,x)\,dx\,dt - \sum_{i\in\mathbb{Z}} \int_0^{\delta t} \int_{i\delta x}^{(i+1)\delta x} \phi_t(0,i\delta x)\rho_{\delta x,\delta t}(t,x)\,dx\,dt \\
&+\delta t \sum_{k\geq 0}\sum_{i\in\mathbb{Z}} \rho_{i-1}^k w_{i-1}^k \int_{(i-1)\delta x}^{i\delta x} \int_{(i-1)\delta x}^{\xi} \phi_{xx}(k\delta t,s)\,ds\,d\xi + \delta x \sum_{k\geq 0}\sum_{i\in\mathbb{Z}} \rho_i^{k+1} \int_{k\delta t}^{(k+1)\delta t} \int_{k\delta t}^{\tau} \phi_{tt}(l,i\delta x)\,dl\,d\tau \\
&-\sum_{i\in\mathbb{Z}} \int_{i\delta x}^{(i+1)\delta x} \int_{i\delta x}^{x} \rho_0(x)\phi_x(0,l)\,dl\,dx \\
&-\sum_{k\geq 0}\sum_{i\in\mathbb{Z}} \int_{k\delta t}^{(k+1)\delta t} \int_{(i-1)\delta x}^{i\delta x} \left(\int_{k\delta t}^{t} \phi_{xt}(\tau,(i-1)\delta x)\,d\tau\right)\rho_{\delta x,\delta t}(t,x)w_{\delta x,\delta t}(t,x)\,dx\,dt \\
&+\sum_{k\geq 0}\sum_{i\in\mathbb{Z}} \int_{(k+1)\delta t}^{(k+2)\delta t} \int_{i\delta x}^{(i+1)\delta x} \left(\phi_t(k\delta t,i\delta x) - \phi_t((k+1)\delta t,i\delta x)\right)\rho_{\delta x,\delta t}(t,x)\,dx\,dt
\end{aligned} \qquad (5.21)$$

Using (5.21) it is straightforward to obtain (3.9) with $P_1$ defined by (3.11). Notice next that definition (3.3) implies that

$$\sum_{k\geq 0}\sum_{i\in\mathbb{Z}} \beta\left(w_i^k\right)\phi_i^k = \sum_{k\geq 0}\sum_{i\in\mathbb{Z}} Q'(\rho_{i+1}^k) \frac{\phi_{i+1}^k - \phi_i^k}{\delta x} \qquad (5.22)$$

and due to (5.20) we obtain that

$$\begin{aligned}
&\sum_{k\geq 0}\sum_{i\in\mathbb{Z}} \int_{k\delta t}^{(k+1)\delta t} \int_{i\delta x}^{(i+1)\delta x} \phi(k\delta t,i\delta x)\beta\left(w_{\delta x,\delta t}(t,x)\right)\,dx\,dt \\
&= \sum_{k\geq 0}\sum_{i\in\mathbb{Z}} \int_{k\delta t}^{(k+1)\delta t} \int_{(i+1)\delta x}^{(i+2)\delta x} \phi_x(k\delta t,i\delta x)Q'\left(\rho_{\delta x,\delta t}(t,x)\right)\,dx\,dt \\
&+ \delta t \sum_{k\geq 0}\sum_{i\in\mathbb{Z}} Q'\left(\rho_{i+1}^k\right) \int_{i\delta x}^{(i+1)\delta x} \int_{i\delta x}^{\xi} \phi_{xx}(k\delta t,s)\,ds\,d\xi
\end{aligned} \qquad (5.23)$$

From (5.23) and expanding the integrals $\int_{k\delta t}^{t} \phi_t(\tau, i\delta x)\,d\tau$ and $\int_{i\delta x}^{x} \phi_x(t,\xi)\,d\xi$ we can obtain the

following equation

$$\sum_{k\geq 0}\sum_{i\in\mathbb{Z}} \int_{k\delta t}^{(k+1)\delta t} \int_{i\delta x}^{(i+1)\delta x} \phi(t,x)\beta\bigl(w_{\delta x,\delta t}(t,x)\bigr)dxdt$$

$$-\sum_{k\geq 0}\sum_{i\in\mathbb{Z}} \int_{k\delta t}^{(k+1)\delta t} \int_{i\delta x}^{(i+1)\delta x} \left(\int_{i\delta x}^{x} \phi_x(t,\xi)d\xi\right)\beta\bigl(w_{\delta x,\delta t}(t,x)\bigr)dxdt$$

$$-\sum_{k\geq 0}\sum_{i\in\mathbb{Z}} \int_{k\delta t}^{(k+1)\delta t} \int_{i\delta x}^{(i+1)\delta x} \left(\int_{k\delta t}^{t} \phi_t(\tau,i\delta x)d\tau\right)\beta\bigl(w_{\delta x,\delta t}(t,x)\bigr)dxdt$$

$$= \sum_{k\geq 0}\sum_{i\in\mathbb{Z}} \int_{k\delta t}^{(k+1)\delta t} \int_{i\delta x}^{(i+1)\delta x} \phi_x(t,i\delta x)Q'\bigl(\rho_{\delta x,\delta t}(t,x)\bigr)dxdt$$

$$-\sum_{k\geq 0}\sum_{i\in\mathbb{Z}} \int_{k\delta t}^{(k+1)\delta t} \int_{(i+1)\delta x}^{(i+2)\delta x} \bigl(\phi_x(t,(i+1)\delta x)-\phi_x(t,i\delta x)\bigr)Q'\bigl(\rho_{\delta x,\delta t}(t,x)\bigr)dxdt$$

$$+\delta t\sum_{k\geq 0}\sum_{i\in\mathbb{Z}} Q'\bigl(\rho_{i+1}^k\bigr) \int_{i\delta x}^{(i+1)\delta x}\int_{i\delta x}^{\xi} \phi_{xx}(k\delta t,s)dsd\xi$$

$$-\sum_{k\geq 0}\sum_{i\in\mathbb{Z}} \int_{k\delta t}^{(k+1)\delta t} \int_{(i+1)\delta x}^{(i+2)\delta x} \left(\int_{k\delta t}^{t} \phi_{xt}(\tau,i\delta x)d\tau\right)Q'\bigl(\rho_{\delta x,\delta t}(t,x)\bigr)dxdt \qquad (5.24)$$

Equation (5.24) implies (3.10) with $P_2$ defined by (3.12). Finally, inequalities (3.13) and (3.14) are a direct consequence of Proposition 3 and Proposition 4. The proof is complete. ◁

**Proof of Proposition 7:** Define

$$V_k = \sum_{i\in\mathbb{Z}} Q(\rho_i^k) \qquad (5.25)$$

$$W_k = \sum_{i\in\mathbb{Z}} \rho_i^k w_i^k \beta\bigl(w_i^k\bigr) \qquad (5.26)$$

and notice that due to definitions (2.19) and (2.41), $V_k \geq 0$ and $W_k \geq 0$ for all $k \geq 0$. Moreover, let

$$c = 1 - \delta t \frac{4\|\rho_0\|_\infty \|h'\|_\infty}{(\delta x)^2} \max_{0\leq \rho \leq \|\rho_0\|_\infty} \bigl(\kappa(\rho)\bigr) > 0 \qquad (5.27)$$

Definitions (5.25), (5.26), (5.27) and inequality (3.4) imply the following estimate

$$V_{k+1} \leq V_k - c\,\delta t\, W_k$$

which recursively gives that

$$V_k \leq V_0 - c\,\delta t \sum_{j=0}^{k-1} W_j \qquad (5.28)$$

Since $V_k \geq 0$, we obtain from (5.28) the following implication

$$\sum_{j=0}^{+\infty} W_j \leq \frac{1}{c \, \delta t} V_0 \Rightarrow \lim_{k \to +\infty} \left( \sum_{i \in \mathbb{Z}} \rho_i^k w_i^k \beta(w_i^k) \right) = 0 \Rightarrow \lim_{k \to +\infty} \left( \rho_i^k w_i^k \beta(w_i^k) \right) = 0 \qquad (5.29)$$

Suppose that $\lim_{k \to +\infty} \left( \rho_i^k |w_i^k| \right) = 0$ is not true for some $i \in \mathbb{Z}$. Then there exists $\varepsilon \in \left( 0, \frac{\min(1,b)}{R} \right)$ and an increasing sequence $\{k_n : n = 1, 2, \ldots\}$ with $\rho_i^{k_n} |w_i^{k_n}| \geq \varepsilon$ for all $n = 1, 2, \ldots$ . It follows that $|w_i^{k_n}| \geq \frac{\varepsilon}{R}$ and $\rho_i^{k_n} \geq \frac{\varepsilon}{\max(1,b)}$ for all $n = 1, 2, \ldots$ . Since $\beta$ is increasing with $\beta(0) = 0$ we get that $w_i^{k_n} \beta(w_i^{k_n}) \geq \min\left( \frac{\varepsilon}{R} \beta\left( \frac{\varepsilon}{R} \right), -\frac{\varepsilon}{R} \beta\left( -\frac{\varepsilon}{R} \right) \right)$ for all $n = 1, 2, \ldots$ . Consequently, we have $\rho_i^{k_n} w_i^{k_n} \beta(w_i^{k_n}) \geq \frac{\varepsilon}{\max(1,b)} \min\left( \frac{\varepsilon}{R} \beta\left( \frac{\varepsilon}{R} \right), -\frac{\varepsilon}{R} \beta\left( -\frac{\varepsilon}{R} \right) \right)$ for all $n = 1, 2, \ldots$, which contradicts the fact that $\lim_{k \to +\infty} \left( \rho_i^k w_i^k \beta(w_i^k) \right) = 0$. Therefore, $\lim_{k \to +\infty} \left( \rho_i^k |w_i^k| \right) = 0$ for some $i \in \mathbb{Z}$.

Suppose that $\lim_{k \to +\infty} \left( \rho_i^k |Q'(\rho_i^k) - Q'(\rho_{i+1}^k)| \right) = 0$ is not true for some $i \in \mathbb{Z}$. Clearly, this cannot happen if $\|\rho_0\|_\infty \leq 1$. Consequently, $\|\rho_0\|_\infty > 1$. Then there exists $\varepsilon > 0$ and an increasing sequence $\{k_n : n = 1, 2, \ldots\}$ with $\rho_i^{k_n} |Q'(\rho_i^{k_n}) - Q'(\rho_{i+1}^{k_n})| \geq \varepsilon$ for all $n = 1, 2, \ldots$ . It follows that $|Q'(\rho_i^{k_n}) - Q'(\rho_{i+1}^{k_n})| \geq \frac{\varepsilon}{R}$ and $\rho_i^{k_n} \geq \frac{\varepsilon}{Q'(\|\rho_0\|_\infty)}$ for all $n = 1, 2, \ldots$ . Since $h$ is increasing with $h(0) = 0$ and since $w_i^k = h\left( \frac{Q'(\rho_i^k) - Q'(\rho_{i+1}^k)}{\delta x} \right)$ we get that $|w_i^{k_n}| \geq \min\left( h\left( \frac{\varepsilon}{R \delta x} \right), -h\left( -\frac{\varepsilon}{R \delta x} \right) \right)$ for all $n = 1, 2, \ldots$ . Consequently, we have $\rho_i^{k_n} |w_i^{k_n}| \geq \frac{\varepsilon}{Q'(\|\rho_0\|_\infty)} \min\left( h\left( \frac{\varepsilon}{R \delta x} \right), -h\left( -\frac{\varepsilon}{R \delta x} \right) \right)$ for all $n = 1, 2, \ldots$, which contradicts the fact that $\lim_{k \to +\infty} \left( \rho_i^k |w_i^k| \right) = 0$. Therefore, $\lim_{k \to +\infty} \left( \rho_i^k |Q'(\rho_i^k) - Q'(\rho_{i+1}^k)| \right) = 0$ for some $i \in \mathbb{Z}$. The proof is complete. ◁